\documentclass[reqno]{amsart}
\usepackage{amsfonts}
\usepackage{color}
\usepackage{diagrams}
\usepackage{amssymb}
\usepackage{amsmath}

\newtheorem{theorem}{Theorem}[section]

\newtheorem{corollary}[theorem]{Corollary}

\newtheorem{definition}[theorem]{Definition}

\newtheorem{lemma}[theorem]{Lemma}

\newtheorem{proposition}[theorem]{Proposition}

\theoremstyle{remark}
\newtheorem{remark}[theorem]{Remark}

\newcommand{\diff}{ {\rm{d}} }

\numberwithin{equation}{section}

\begin{document}
\title[Dynamic Boundary Conditions]{Hyperbolic Relaxation of Reaction
Diffusion Equations with Dynamic Boundary Conditions}
\author[C. G. Gal and J. L. Shomberg]{Ciprian G. Gal and Joseph L. Shomberg}
\address{Ciprian G. Gal, Department of Mathematics, Florida International
University, Miami, FL 33199, USA}
\email{cgal@fiu.edu}
\address{Joseph L. Shomberg, Department of Mathematics and Computer Science,
Providence College, Providence, RI 02918, USA}
\email{jshomber@providence.edu}
\date{today}

\begin{abstract}
Under consideration is the hyperbolic relaxation of a semilinear
reaction-diffusion equation,
\begin{equation*}
\varepsilon u_{tt}+u_{t}-\Delta u+f(u)=0,
\end{equation*}%
on a bounded domain $\Omega \subset \mathbb{R}^{3}$, with $\varepsilon \in
(0,1]$ and the prescribed dynamic condition,
\begin{equation*}
\partial _{\mathbf{n}}u+u+u_{t}=0,
\end{equation*}%
on the boundary $\Gamma :=\partial \Omega $. We also consider the limit
parabolic problem ($\varepsilon =0$) with the same dynamic boundary
condition. Each problem is well-posed in a suitable phase space where the
global weak solutions generate a Lipschitz continuous semiflow which admits
a bounded absorbing set. Because of the nature of the boundary condition,
fractional powers of the Laplace operator are not well-defined. The
precompactness property required by the hyperbolic semiflows for the
existence of the global attractors is gained through the approach of \cite%
{Pata&Zelik06}. In this case, the optimal regularity for the global
attractors is also readily established. In the parabolic setting, the
regularity of the global attractor is necessary for the semicontinuity
result. After fitting both problems into a common framework, a proof of the
upper-semicontinuity of the family of global attractors is given at $%
\varepsilon =0$. Finally, we also establish the existence of a family of
exponential attractors.
\end{abstract}

\maketitle
\tableofcontents

\section{Introduction}


Let $\Omega $ be a bounded domain in $\mathbb{R}^{3}$ with boundary $\Gamma
:=\partial \Omega $ of class $C^{2}$. We consider the \emph{hyperbolic}
relaxation of a semilinear reaction diffusion equation,
\begin{equation}
\varepsilon u_{tt}+u_{t}-\Delta u+f(u)=0,  \label{eq:pde-h-1}
\end{equation}%
in$~~(0,\infty )\times \Omega ,$ where $\varepsilon \in \lbrack 0,1]$. The
equation is endowed with the dynamic boundary condition,
\begin{equation}
\partial _{\mathbf{n}}u+u+u_{t}=0,  \label{eq:pde-h-2}
\end{equation}%
on$~~(0,\infty )\times \Gamma ,$ and with the initial conditions,
\begin{equation}
u(0,x)=u_{0}(x),~~u_{t}(0,x)=u_{1}(x)~~\text{in}~\Omega .  \label{eq:pde-h-3}
\end{equation}

For the nonlinear term $f$, we assume, $f\in C^{2}(\mathbb{R})$ and that
there is a constant, $\ell \geq 0$, such that, for all $s\in \mathbb{R}$,
the following growth and sign conditions are satisfied,
\begin{equation}
|f^{\prime \prime }(s)|\leq \ell (1+|s|),  \label{eq:assumption-f-1}
\end{equation}%
and
\begin{equation}
\liminf_{|s|\rightarrow \infty }\frac{f(s)}{s}>-\lambda ,
\label{eq:assumption-f-2}
\end{equation}%
where $\lambda >0$ is the best Sobolev/Poincar\'{e} type constant
\begin{equation}
\lambda \int_{\Omega }u^{2}\mathrm{d}x\leq \int_{\Omega }|\nabla u|^{2}%
\mathrm{d}x+\int_{\Gamma }u^{2}\mathrm{d}S.  \label{eq:Poincare-type}
\end{equation}%
Finally, assume that there is $\vartheta >0$ such that for all $s\in \mathbb{%
R}$,
\begin{equation}
f^{\prime }(s)\geq -\vartheta .  \label{eq:assumption-f-3}
\end{equation}%
Notice that the derivative $f=F^{^{\prime }}$\ of the double-well potential,
$F(u)=\frac{1}{4}u^{4}-ku^{2}$, $k>0$, satisfies assumptions (\ref%
{eq:assumption-f-1}), (\ref{eq:assumption-f-2}) and (\ref{eq:assumption-f-3}%
). The first two assumptions made here on the nonlinear term, (\ref%
{eq:assumption-f-1}) and (\ref{eq:assumption-f-2}), are the same assumptions
made on the nonlinear term in \cite{CEL02}, \cite{MPZ07} and \cite%
{Wu&Zheng06}, for example (\cite{MPZ07} additionally assumes $f(0)=0$). The
third assumption (\ref{eq:assumption-f-3}) appears in \cite{CGG11}, \cite%
{Frigeri10}, \cite{GM12} and \cite{Pata&Zelik06}; the bound is utilized to
obtain the precompactness property for the semiflow associated with
evolution equations when dynamic boundary conditions present a difficulty
(e.g., here, fractional powers of the Laplace operator subject to (\ref%
{eq:pde-h-2}) are undefined). It is worth mentioning that (\ref%
{eq:assumption-f-2}) can be also replaced by a less general (but still
widely used in the literature) condition%
\begin{equation*}
\liminf_{|s|\rightarrow \infty }f^{^{\prime }}(s)\geq -\lambda
\end{equation*}%
in which case, (\ref{eq:assumption-f-3}) is automatically satisfied.
Furthermore, assumption (\ref{eq:assumption-f-1}) implies that the growth
condition for $f$ is the critical case since $\Omega \subset \mathbb{R}^{3}$%
. Such assumptions are common when one is investigating the existence of a
global attractor or the existence of an exponential attractor for a partial
differential equation of evolution.

Of course, when (\ref{eq:pde-h-1}) is equipped with Dirichlet, Neumann or
periodic boundary conditions, (\ref{eq:Poincare-type}) simplifies. Moreover,
if (\ref{eq:pde-h-1}) is equipped with a Robin boundary condition, then an
estimate like (\ref{eq:Poincare-type}) holds, but $\lambda $ possesses an
explicit description as the first eigenvalue of the Laplacian with respect
to the Robin boundary condition. The relation between the dynamic condition (%
\ref{eq:pde-h-2}) with the acoustic boundary condition is discussed below.
The hyperbolic equation (\ref{eq:pde-h-1}) is a well-known nonlinear wave
equation motivated from (relativistic) quantum mechanics (cf., e.g \cite%
{Babin&Vishik92,Chepyzhov&Vishik02,Ladyzhenskaya91,Temam88}). However, as
mentioned, most sources study the asymptotic behavior of (\ref{eq:pde-h-1})
with a static boundary condition such as Dirichlet, Neumann, periodic or
Robin. One of the goals of this paper is to extend some results concerning
the asymptotic behavior of (\ref{eq:pde-h-1}), now with the dynamic boundary
condition (\ref{eq:pde-h-2}). The corresponding linear case for (\ref%
{eq:pde-h-1})-(\ref{eq:pde-h-3}) is treated in \cite{PRB04}. The existence
of the global attractor for a linear damped wave equation with a nonlinear
dynamic boundary condition is considered in \cite{Z}. More general systems,
with supercritical nonlinear sources on both the interior and the boundary,
are considered in \cite{B1, BL1, BL2, BL3, BL4}. These contributions mainly
devote their attention to issues like, Hadamard local wellposedness, global
existence, blow-up and non-existence theorems, as well as estimates on the
uniform energy dissipation rates for the appropriate classes of solutions.
We also refer the reader to \cite{BRT} for a unified overview of these
results.

Our main goal is to compare the hyperbolic relaxation problem (\ref%
{eq:pde-h-1})-(\ref{eq:pde-h-3}) with that of the limit parabolic equation
where, for $\varepsilon =0$, we have the reaction-diffusion equation,
\begin{equation}
u_{t}-\Delta u+f(u)=0,  \label{eq:pde-p-1}
\end{equation}%
in$~~(0,\infty )\times \Omega ,$ with the dynamic boundary condition,
\begin{equation}
\partial _{\mathbf{n}}u+u+u_{t}=0,  \label{eq:pde-p-2}
\end{equation}%
on$~~(0,\infty )\times \Gamma $, and the initial conditions%
\begin{equation}
u(0,x)=u_{0}(x)~~\text{in}~~\Omega ,\text{ }u(0,x)=\gamma _{0}(x)~~\text{on}%
~\Gamma .  \label{eq:pde-p-3}
\end{equation}%
For the sake of simplicity, we shall restrict our attention only to linear
boundary conditions of the form (\ref{eq:pde-p-2}) even though our framework
can easily allow for a complete treatment of nonlinear dynamic boundary
conditions (see Remark \ref{nonlinear_bc}; cf. also \cite{CGG11}, \cite%
{Frigeri10}, \cite{GM12}).

Because of its importance in the physical sciences and the development of
mathematical physics, the reaction-diffusion equation (\ref{eq:pde-p-1}) and
its asymptotic behavior are well-known to the literature. Many of the books
referenced above contain a treatment on the parabolic semilinear
reaction-diffusion equation (\ref{eq:pde-p-1}) with \emph{static} boundary
conditions. In particular, the Chaffee-Infante reaction-diffusion equation,
with $f(u)=u^3-ku$, $k>0$ and Dirichlet boundary conditions can be found in
\cite[Section 11.5]{Robinson01}. A discussion on the structure of the
associated global attractor can also be found there. Additionally, the
Chaffee-Infante equation and its hyperbolic relaxation, again with Dirichlet
boundary conditions, are discussed in \cite[Chapters 3-5]{Milani&Koksch05}.

Recently there has been a great amount of research taking place in the area
of partial differential equations of evolution type, subject to dynamic
boundary conditions. Boundary conditions of the form (\ref{eq:pde-p-2})
arise for many known equations of mathematical physics. This can especially
be seen by the many applications given to heat control problems,
phase-transition phenomena, Stefan problems, some models in climatology, and
many others. Without being too exhaustive we refer the reader to \cite%
{Gal12-2,Gal12,Gal&Warma10} for more details about the system (\ref%
{eq:pde-p-1})-(\ref{eq:pde-p-3})\ and a more complete list of references. A
version of equation (\ref{eq:pde-h-2}), but with nonlinear dissipation on
the boundary, already appears in the literature, we refer to \cite%
{CEL02,CEL04-2,CEL04}. There, the authors are able to show the existence of
a global attractor without the presence of the weak interior damping term $%
u_{t},$ by assuming that $f$ is \emph{subcritical}. One motivation for
considering a boundary condition like (\ref{eq:pde-h-2}) comes from
mechanical considerations: there is frictional damping on the boundary $%
\Gamma $ that is linearly proportional to the velocity $u_{t}$. In \cite%
{Wu&Zheng06}, the convergence, as time goes to infinity, of unique global
strong solutions of (\ref{eq:pde-h-1})-(\ref{eq:pde-h-3})\ to a single
equilibrium is established provided that $f$ is also real analytic. Note
that the set of equilibria for (\ref{eq:pde-h-1})-(\ref{eq:pde-h-3}) may
form a continuum so that, in general, guaranteeing this convergence is a
highly nontrivial matter. The second motivation comes from thermodynamics.
Suppose that we want to consider heat flow in a metal. The standard
derivation of the heat equation is always based on the idea that
\textquotedblleft heat in equals heat out\textquotedblright\ over a region $%
\overline{\Omega }$. But the classical approach ignores the contribution of
heat sources located on the boundary $\Gamma $, by taking into account only
heat sources/sinks which are present inside the region (in our case, $%
-f\left( u\right) $ is treated as a source within $\Omega $). A new
derivation of the heat equation in the presence of heat sources/sinks
located at $\Gamma ,$ assuming the Fourier law of cooling states (i.e., the
heat flux $\overrightarrow{q}$ is directly proportional to the temperature
gradient, $\overrightarrow{q}=-\nabla u$) was given in \cite{Gold}, and it
has lead to the precise formulation of the system in (\ref{eq:pde-p-1})-(\ref%
{eq:pde-p-3}). However, the derivation in \cite{Gold} suffers from an
important drawback which cannot be ignored: initial perturbations in (\ref%
{eq:pde-p-1}) propagate with infinite speed. This means that the presence of
a heat source located at $\Gamma $ is instantaneously felt by all observers
in $\Omega $, no matter how far away from $\Gamma $ they happen to be. This
behavior can be traced to the \textquotedblleft parabolic\textquotedblright\
character of Fourier's law. Thus, in many relevant phenomena the system (\ref%
{eq:pde-p-1})-(\ref{eq:pde-p-3}) can become a bad approximation (see, e.g.,
\cite{APR}, \cite{HP}, for many examples). In order to overcome these
problems, a generalization of the standard Fourier law must be considered,
leading to a new formulation for which the heat flux $\overrightarrow{q}$
obeys the so-called Maxwell--Cattaneo heat conduction law:%
\begin{equation}
\varepsilon \partial _{t}\overrightarrow{q}+\overrightarrow{q}=-\nabla u,
\label{MClaw}
\end{equation}%
in $\left( 0,\infty \right) \times \Omega .$ Note that the Fourier law is
obtained from (\ref{MClaw}) when $\varepsilon =0$. This expression for the
heat flux $\overrightarrow{q}$ leads to the hyperbolic equation (\ref%
{eq:pde-h-1}), which entails that $u$ propagates at finite speed. It is also
worth mentioning that one can write (\ref{MClaw}) in the equivalent form of%
\begin{equation}
\overrightarrow{q}\left( t,x\right) =-\int_{0}^{\infty }\Theta _{\varepsilon
}\left( t-s\right) \nabla u\left( s,x\right) \mathrm{d}s,\text{ }\Theta
_{\varepsilon }\left( t\right) :=\frac{1}{\varepsilon }e^{-\frac{t}{%
\varepsilon }}.  \label{MClaw2}
\end{equation}%
This points to a situation in which the (past) thermal memory of the
material plays a role, but its relevance goes down quickly as we move to the
past. Finally, it may be worth mentioning that the form of flux $%
\overrightarrow{q}$ assumed in (\ref{MClaw2}), in which $\Theta
_{\varepsilon }$ is assumed to be a \emph{generic} memory kernel, also
yields the following problem:%
\begin{equation}
u_{t}=\int_{0}^{\infty }\Theta _{\varepsilon }\left( t-s\right) \left(
\Delta u\left( s\right) -f\left( u\left( s\right) \right) \right) \mathrm{d}%
s.  \label{mem-eq}
\end{equation}%
In this case, $\Theta _{\varepsilon }\left( s\right) =\varepsilon
^{-1}\Theta \left( s/\varepsilon \right) $ and $\Theta :\left( 0,\infty
\right) \rightarrow \left( 0,\infty \right) $ is a given (smooth) summable
and convex (hence decreasing) relaxation kernel. A complete treatment of
equation (\ref{mem-eq}), endowed with the dynamic boundary condition (\ref%
{eq:pde-p-2}), will be the subject of further investigation in the future.

It may also be interesting to note that the dynamic boundary condition given
in (\ref{eq:pde-h-2}) can be recovered, in some sense, from the linear
acoustic boundary condition,
\begin{equation}
\left\{
\begin{array}{ll}
m\delta _{tt}+\delta _{t}+\delta =-u_{t} & \text{on}~~(0,T)\times \Gamma \\
\partial _{\mathbf{n}}u=\delta _{t} & \text{on}~~(0,T)\times \Gamma .%
\end{array}%
\right.  \label{eq:acoustic}
\end{equation}%
Here the unknown $\delta =\delta (t,x)$ represents the \emph{inward}
\textquotedblleft displacement\textquotedblright\ of the boundary $\Gamma $
reacting to a pressure described by $-u_{t}$. The first equation (\ref%
{eq:acoustic})$_{1}$ describes the spring-like effect in which $\Gamma $
(and $\delta $) interacts with $-u_{t}$, and the second equation (\ref%
{eq:acoustic})$_{2}$ is the continuity condition: velocity of the boundary
displacement $\delta $ agrees with the normal derivative of $u$. Together, (%
\ref{eq:acoustic}) describes $\Gamma $ as a locally reactive surface. The
term $m=m(x)$ represents mass, so in a mass\emph{less} system, the inertial
term disappears. In the case when $\delta $ can be modelled by $u$ near the
boundary; i.e., if $\delta \sim u$ near $\Gamma $, then we arrive at the
boundary condition described by (\ref{eq:pde-h-2}). In applications, the
unknown $u$ may be taken as a velocity potential of some fluid or gas in $%
\Omega $ that was disturbed from its equilibrium. The acoustic boundary
condition was rigorously described by Beale and Rosencrans in \cite{Beale76}
and \cite{Beale&Rosencrans74}. Various recent sources investigate the wave
equation equipped with acoustic boundary conditions, \cite%
{CFL01,GGG03,Mugnolo10,Vicente09}. However, more recently, it has been
introduced as a dynamic boundary condition for problems that study the
asymptotic behavior of weakly damped wave equations, see \cite{Frigeri10}
and \cite{Shomberg&Frigeri12}.

The aim of this paper is to extend the asymptotic results for dissipative
wave equations (\ref{eq:pde-h-1}) and reaction-diffusion equations (\ref%
{eq:pde-p-1})\ with the dynamic boundary condition (\ref{eq:pde-h-2}), in
terms of a perturbation problem, and ultimately discuss the continuity of
the attracting sets generated by these problems. Due to the nature of the
boundary condition imposed for the model problem (\ref{eq:pde-h-1}), we are
unable to prove the existence of global attractors for the hyperbolic
relaxation problem through the compactness argument which is typical for
damped wave equations with static boundary conditions, such as Dirichlet,
Neumann, periodic, or Robin boundary conditions (cf. e.g. \cite%
{Milani&Koksch05,Temam88,Zheng04}). The problem arises from our lack to
define fractional powers of the Laplacian with respect to the boundary
condition (\ref{eq:pde-h-2}). This situation takes place because of the
permanence of the $u_{t}$ term on $\Gamma $, which in turn means the
\textquotedblleft Laplacian\textquotedblright\ is not self-adjoint. Thus,
for example, the model problem does not enjoy an explicit Poincar\'{e}
inequality found with a Fourier series, nor the existence of a local weak
solution found with a typical Galerkin basis. Local solutions will be sought
with semigroup methods that rely on monotone operators techniques as in \cite%
{CEL02}. Then estimates are applied to extend the local solutions to global
ones and the existence of an absorbing set is determined. For the hyperbolic
relaxation problem (\ref{eq:pde-h-1})-(\ref{eq:pde-h-3}), we obtain the
relatively compact part in the decomposition of the solution by following
the approach in \cite{Pata&Zelik06}.

The main novelties of the present paper with respect to previous results on
the damped wave equation (\ref{eq:pde-h-1}) are the following:

\begin{itemize}
\item We extend the results on the existence of global attractors $\left\{
\mathcal{A}_{\varepsilon }\right\} _{\varepsilon \in (0,1]}$ for the damped
wave equation (\ref{eq:pde-h-1}) with a critical nonlinearity and a
``dynamic'' boundary condition instead of the usual Dirichlet boundary
condition (see, e.g., \cite{Hale&Raugel88}, \cite{Hale88}). This is achieved
through the decomposition method exploited in \cite{Pata&Zelik06} which
allows us to establish that $\mathcal{A}_{\varepsilon }$ has also \emph{%
optimal} regularity (see Theorem \ref{t:global-attractors-h}).

\item We show that a certain family $\left\{ \widetilde{\mathcal{A}}%
_{\varepsilon }\right\} _{\varepsilon \in \left[ 0,1\right] }$\ of compact
sets, which is topologically conjugated to $\left\{ \mathcal{A}_{\varepsilon
}\right\} _{\varepsilon \in \left[ 0,1\right] }$ in a precise way, is also
upper-semicontinuous as $\varepsilon $ goes to zero. Roughly speaking, we
show that these sets $\widetilde{\mathcal{A}}_{\varepsilon }$ converge to
the \textquotedblleft lifted\textquotedblright\ global attractor $\widetilde{%
\mathcal{A}}_{0}$ associated with the parabolic problem. The argument
utilizes the sequential characterization of the global attractor (cf., e.g.
\cite[Proposition 2.15]{Milani&Koksch05}). The main difficulty comes from
the fact that the phase spaces for the perturbed and unperturbed equations
are not the same; indeed, solutions of the hyperbolic problem are defined
for $\left( u_{0},u_{1}\right) \in H^{s+1}\left( \Omega \right) \times
H^{s}\left( \Omega \right) $, $s\in \left\{ 0,1\right\} $, while solutions
of the parabolic problem make sense only in spaces like \thinspace $%
L^{2}\left( \Omega \right) \times L^{2}\left( \Gamma \right) $ and $%
H^{s+1}\left( \Omega \right) \times H^{s+1/2}\left( \Gamma \right) ,$
respectively (see (\ref{eq:pde-p-3})). Thus, previous constructions obtained
for parabolic equations with Dirichlet boundary conditions cannot be applied
and have to be adapted.

\item We prove the existence of a family of exponential attractors $\left\{
\mathcal{M}_{\varepsilon }\right\} ,$ $\varepsilon \in (0,1]$, which entails
that $\mathcal{A}_{\varepsilon }$ is also finite dimensional even in the
critical case. We recall that the same result was shown in \cite{CEL02} for
the wave equation (i.e., (\ref{eq:pde-h-1}) without any damping in $\Omega $%
) subject to the boundary condition (\ref{eq:pde-p-2}) only in the
subcritical case. Unfortunately, we are unable to show that this dimension
is uniform with respect to $\varepsilon >0$ as $\varepsilon $ goes to zero.
Some other open questions are formulated at the end of the article.
\end{itemize}

The article is organized as follows. The limit ($\varepsilon =0$)
reaction-diffusion problem is discussed in Section \ref{s:parabolic}. The
section is mostly devoted to citing the already known main results of the
parabolic problem: the existence and uniqueness of global solutions in an
appropriate phase space (see Theorem \ref{t:parabolic-weak-solutions}), the
definition of the (Lipschitz) semiflow, the existence and regularity of the
global attractor (see Theorem \ref{t:global-attractor-and-regularity-p}).
Section \ref{s:hyperbolic} contains our treatment of the hyperbolic
relaxation problem, for all $\varepsilon \in (0,1]$. We discuss the
existence and uniqueness of solutions defined for all positive times in
Section \ref{well-posed} (see Theorem \ref{t:hyperbolic-weak-solutions}).
The solutions generate a semiflow on the phase space, and thanks to the
continuous dependence estimate, we know that the semiflow is locally
Lipschitz continuous. The existence of a bounded absorbing set is also shown
(see Lemma \ref{t:bounded-absorbing-set-h}). The global attractor and its
properties are established in Section \ref{global-hyper}, while the
upper-semicontinuous result is established in Section \ref{s:continuity}.
The existence of exponential attractors for the hyperbolic problem is
presented in Section \ref{s:exponential-attractors}. The statement of a Gr%
\"{o}nwall-type inequality, used frequently in the estimates, is included in
the Appendix.

\section{The limit parabolic problem}

\label{s:parabolic}

In this short section, we recall some results for the limit parabolic
problem (\ref{eq:pde-p-1})-(\ref{eq:pde-p-3}), i.e., (\ref{eq:pde-h-1})-(\ref%
{eq:pde-h-3}) with $\varepsilon =0$. Unlike the hyperbolic problem, a full
general treatment of the limit parabolic problem with dynamic boundary
conditions already appears in the literature (cf., e.g., \cite%
{Gal12-2,Gal12,Gal&Warma10, Mey10} and references there in); in particular,
this section will summarize some of the main results from \cite{Gal12}. It
should be noted for the interest of the reader that all formal calculations
made with the weak solutions of the parabolic problem can be rigorously
justified using the Galerkin discretization scheme that appears, for
instance, in \cite[Theorem 2.6]{Gal&Warma10}. Indeed, it is through the use
of the Galerkin approximations that the existence of weak solutions for the
parabolic problem is shown. The solution operator associated with the
parabolic problem generates a locally Lipschitz continuous semiflow on the
appropriate phase space. We also know that this semiflow admits a connected
global attractor that is bounded in a more \emph{regular} phase space. It
follows that solutions, when restricted to the global attractor, are in fact
strong solutions, exhibiting further regularity that will become \emph{%
essential} when we later consider the continuity properties of the family of
global attractors produced by the hyperbolic relaxation problem ($%
\varepsilon >0$) and the limit parabolic problem ($\varepsilon =0$).

We need to introduce some notations and definitions. From now on, we denote
by $\Vert \cdot \Vert $, $\Vert \cdot \Vert _{k}$, the norms in $%
L^{2}(\Omega )$, $H^{k}(\Omega )$, respectively. We use the notation $%
\langle \cdot ,\cdot \rangle $ and $\langle \cdot ,\cdot \rangle _{k}$ to
denote the products on $L^{2}(\Omega )$ and $H^{k}(\Omega )$, respectively.
For the boundary terms, $\Vert \cdot \Vert _{L^{2}(\Gamma )}$ and $\langle
\cdot ,\cdot \rangle _{L^{2}(\Gamma )}$ denote the norm and, respectively,
product on $L^{2}(\Gamma )$. We will require the norm in $H^{k}(\Gamma )$,
to be denoted by $\Vert \cdot \Vert _{H^{k}(\Gamma )}$, where $k\geq 1$. The
$L^{p}(\Omega )$ norm, $p\in (0,\infty ]$, is denoted $|\cdot |_{p}$. The
dual pairing between $H^{1}(\Omega )$ and its dual $(H^{1}(\Omega ))^{\ast }$
is denoted by $(u,v)$. We denote the measure of the domain $\Omega $ by $%
|\Omega |$. In many calculations, functional notation indicating dependence
on the variable $t$ is dropped; for example, we will write $u$ in place of $%
u(t)$. Throughout the paper, $C\geq 0$ will denote a \emph{generic}
constant, while $Q:\mathbb{R}_{+}\rightarrow \mathbb{R}_{+}$ will denote a
\emph{generic} increasing function. All these quantities, unless explicitly
stated, are \emph{independent} of $\varepsilon .$ Further dependencies of
these quantities will be specified on occurrence.

The following inequalities are straight forward consequences of the Poincar%
\'{e} type inequality (\ref{eq:Poincare-type}) and assumptions (\ref%
{eq:assumption-f-2}) and (\ref{eq:assumption-f-3}). From (\ref%
{eq:assumption-f-2}) it follows that, for some constants $\mu \in (0,\lambda
]$ and $c_{1}=c_{1}(f,|\Omega |)\geq 0$, and for all $\xi \in H^{1}(\Omega )$%
,%
\begin{align}
\left\langle f\left( \xi \right) ,\xi \right\rangle & \geq -\left( \lambda
-\mu \right) \left\Vert \xi \right\Vert ^{2}-c_{1}
\label{eq:consequence-f-1} \\
& \geq -\frac{\left( \lambda -\mu \right) }{\lambda }\left( \left\Vert
\nabla \xi \right\Vert ^{2}+\left\Vert \xi \right\Vert _{L^{2}\left( \Gamma
\right) }^{2}\right) -c_{1}.  \notag
\end{align}%
Let $F(s)=\int_{0}^{s}f(\sigma )\mathrm{d}\sigma $. For some constant $%
c_{2}=c_{2}(f,|\Omega |)\geq 0$, and for all $\xi \in H^{1}(\Omega )$,
\begin{equation}
\begin{aligned} \int_\Omega F(\xi)\mathrm{d} x & \geq
-\frac{\lambda-\mu}{2}\|\xi\|^2 - c_2 \\ & \geq
-\frac{\lambda-\mu}{2\lambda}\|\xi\|^2_1 - c_2. \label{eq:consequence-F-1}
\end{aligned}
\end{equation}%
See \cite[page 1913]{CEL02} for an explicit proof of (\ref%
{eq:consequence-F-1}). The proof of (\ref{eq:consequence-f-1}) is similar.
Finally, using (\ref{eq:assumption-f-3}) and integration by parts on $%
F(s)=\int_{0}^{s}f(\sigma )\mathrm{d}\sigma $, we have the upper-bound
\begin{equation}
\begin{aligned} \int_\Omega F(\xi)\mathrm{d} x & \leq \langle f(\xi),\xi
\rangle + \frac{\vartheta}{2}\|\xi\|^2 \\ & \leq \langle f(\xi),\xi \rangle
+ \frac{\vartheta}{2\lambda}\|\xi\|^2_1. \label{eq:consequence-F-2}
\end{aligned}
\end{equation}

The natural energy phase space for the limit parabolic problem (\ref%
{eq:pde-p-1})-(\ref{eq:pde-p-3}) is the space
\begin{equation*}
Y=L^{2}(\Omega )\times L^{2}(\Gamma ),
\end{equation*}%
which is Hilbert when equipped with the norm whose square is given by, for
all $\zeta =(u,\gamma )\in Y$,
\begin{equation*}
\Vert \zeta \Vert _{Y}^{2}:=\Vert u\Vert ^{2}+\Vert \gamma \Vert
_{L^{2}(\Gamma )}^{2}.
\end{equation*}%
It is well-known that the Dirichlet trace map $\mathrm{tr_{D}}:C^{\infty
}\left( \overline{\Omega }\right) \rightarrow C^{\infty }\left( \Gamma
\right) ,$ defined by $\mathrm{tr_{D}}\left( u\right) =u_{\mid \Gamma }$
extends to a linear continuous operator $\mathrm{tr_{D}}:H^{r}\left( \Omega
\right) \rightarrow H^{r-1/2}\left( \Gamma \right) ,$ for all $r>1/2$, which
is onto for $1/2<r<3/2.$ This map also possesses a bounded right inverse $%
\mathrm{tr_{D}^{-1}}:H^{r-1/2}\left( \Gamma \right) \rightarrow H^{r}\left(
\Omega \right) $ such that $\mathrm{tr_{D}}\left( \mathrm{tr_{D}^{-1}}\psi
\right) =\psi ,$ for any $\psi \in H^{r-1/2}\left( \Gamma \right) $.
Identifying each function $\psi \in C\left( \overline{\Omega }\right) $ with
the vector $V=\left( \psi ,\mathrm{tr_{D}}\left( \psi \right) \right) \in
C\left( \overline{\Omega }\right) \times C\left( \Gamma \right) $, it
follows that $C\left( \overline{\Omega }\right) $ is a dense subspace of $%
Y=L^{2}\left( \Omega \right) \times L^{2}\left( \Gamma \right) $ (see, e.g.,
\cite[Lemma 2.1]{DR}). Also, we introduce the subspaces of $H^{r}\left(
\Omega \right) \times H^{r-1/2}\left( \Gamma \right) $, for every $r>1/2$,%
\begin{equation*}
\mathcal{V}^{r}:=\left\{ \left( u,\gamma \right) \in H^{r}\left( \Omega
\right) \times H^{r-1/2}\left( \Gamma \right) :\gamma =\mathrm{tr_{D}}\left(
u\right) \right\} ,
\end{equation*}%
and note that we have the following dense and compact embeddings $\mathcal{V}%
^{r_{1}}\hookrightarrow \mathcal{V}^{r_{2}},$ for any $r_{1}>r_{2}>1/2$. The
linear subspace $\mathcal{V}^{r}$ is densely and compactly embedded into $Y,$
for any $r>1/2$. We emphasize that $\mathcal{V}^{r}$ is not a product space
and that, due to the boundedness of the trace operator $\mathrm{tr_{D}}$,
the space $\mathcal{V}^{r}$ is topologically isomorphic to $H^{r}\left(
\Omega \right) $ in the obvious way. Thus, we can identify each $u\in
H^{r}\left( \Omega \right) $ with a pair $\left( u,\mathrm{tr_{D}}\left(
u\right) \right) \in \mathcal{V}^{r}$. Finally, note that both spaces $%
H^{r}\left( \Omega \right) $ and $\mathcal{V}^{r}$ are normed spaces with
equivalent norms.

The following definition of weak solution to problem (\ref{eq:pde-p-1})-(\ref%
{eq:pde-p-3}) is taken from \cite{Gal&Warma10} (see, e.g., \cite[Definition
2.1]{Gal12-2}, for the more general case).

\begin{definition}
\label{weak-par}Let $T>0$ and $(u_{0},\gamma _{0})\in Y=L^{2}(\Omega )\times
L^{2}(\Gamma )$. The pair $\zeta (t)=(u(t),\gamma (t))$ is said to be a
(global) \emph{weak solution} of (\ref{eq:pde-p-1})-(\ref{eq:pde-p-3}) on $%
[0,T]$ if, for almost all $t\in (0,T]$, $\gamma (t)=u_{\mid \Gamma }(t),$
and $\zeta $ fulfills%
\begin{eqnarray*}
\zeta &\in &C\left( \left[ 0,T\right] ;Y\right) \cap L^{2}\left( 0,T;%
\mathcal{V}^{1}\right) , \\
\partial _{t}\zeta &\in &L^{2}(0,T;\left( \mathcal{V}^{1}\right) ^{\ast }),%
\text{ }u\in H_{\mathrm{loc}}^{1}\left( (0,T];L^{2}\left( \Omega \right)
\right) , \\
\gamma &\in &H_{\mathrm{loc}}^{1}\left( (0,T];L^{2}\left( \Gamma \right)
\right) ,
\end{eqnarray*}%
such that the following identity holds, for almost all $t\in \lbrack 0,T]$
and for all $\xi =(\chi ,\psi )\in \mathcal{V}^{1}$,
\begin{equation}
\left( \partial _{t}\zeta ,\xi \right) _{\left( \mathcal{V}^{1}\right)
^{\ast },\mathcal{V}^{1}}+\langle \nabla u,\nabla \chi \rangle +\langle
f(u),\chi \rangle +\left\langle u,\psi \right\rangle _{L^{2}\left( \Gamma
\right) }=0.  \label{weakf-par}
\end{equation}%
Moreover,%
\begin{equation*}
\zeta (0)=(u_{0},\gamma _{0})=:\zeta _{0}\text{ a.e. in }Y.
\end{equation*}%
The map $\zeta =(u,\gamma )$ is a weak solution on $[0,\infty )$ (i.e. \emph{%
a global weak solution}) if it is a weak solution on $[0,T]$, for all $T>0$.
\end{definition}

\begin{remark}
It is important to observe, for the weak solutions of Definition \ref%
{weak-par}, that $\gamma _{0}=u_{\mid \Gamma }\left( 0\right) $ need not be
the trace of $u_{0}=u_{\mid \Omega }\left( 0\right) $ at the boundary, and
so in this context the boundary equation (\ref{eq:pde-p-2}) is interpreted
as an additional parabolic equation, now acting on the boundary $\Gamma $.
However, the weak solution does fulfill $\gamma (t)=\mathrm{tr_{D}}u(t)$,
for almost all $t>0.$
\end{remark}

The existence part of the (global) weak solutions is from \cite[Theorem 2.6]%
{Gal&Warma10}, and the continuous dependence with respect to the initial
data $\zeta _{0}$, local Lipschitz continuity on $Y$, uniformly in $t$ on
compact intervals, and the uniqueness of the weak solutions follow from \cite%
[Proposition 2.8]{Gal12-2} (cf. also \cite[Lemma 2.7]{Gal&Warma10}).

\begin{theorem}
\label{t:parabolic-weak-solutions} Assume (\ref{eq:assumption-f-1}), (\ref%
{eq:assumption-f-2}) and (\ref{eq:assumption-f-3}) hold. For each $\zeta
_{0}=(u_{0},\gamma _{0})\in Y$, there exists a unique global weak solution
in the sense of Definition \ref{weak-par}. Moreover, the following estimate
holds, for all $t\geq 0$,
\begin{equation}
\Vert \zeta (t)\Vert _{Y}^{2}+\int_{t}^{t+1}\Vert \zeta (s)\Vert _{\mathcal{V%
}^{1}}^{2}\mathrm{d}s\leq C\Vert \zeta _{0}\Vert _{Y}^{2}e^{-\rho t}+C,
\label{dissip-par}
\end{equation}%
for some positive constants $\rho ,C>0$. Furthermore, let $\zeta
(t)=(u(t),\gamma (t))$ and $\theta (t)=(\chi (t),\psi (t))$ denote the
corresponding weak solutions with initial data $\zeta _{0}=(u_{0},\gamma
_{0})$ and $\theta _{0}=(\chi _{0},\psi _{0})$, respectively. Then, for all $%
t\geq 0$,
\begin{equation}
\Vert \zeta (t)-\theta (t)\Vert _{Y}\leq Ce^{\theta t}\Vert \zeta
_{0}-\theta _{0}\Vert _{Y},  \label{eq:continuous-dependence-p}
\end{equation}%
where $C=C\left( R\right) >0$ is such that $\left\Vert \zeta _{0}\right\Vert
_{Y}\leq R,\left\Vert \theta _{0}\right\Vert _{Y}\leq R.$
\end{theorem}

\begin{proof}
Since the proofs in \cite{Gal&Warma10}, \cite{Gal12-2}, \cite{Gal12} involve
quite different assumptions on the nonlinearity other than the ones in the
statement of the theorem, we will sketch a short proof of (\ref{dissip-par}%
). This is the main estimate on which the proof for the existence of a weak
solution is based on (of course, (\ref{dissip-par}) can be rigorously
justified using a suitable Galerkin discretization scheme). To this end,
testing (\ref{weakf-par}) with $\zeta $, and appealing to (\ref%
{eq:consequence-f-1}), we deduce the following inequality%
\begin{equation}
\frac{1}{2}\frac{\mathrm{d}}{\mathrm{d}t}\left\Vert \zeta \left( t\right)
\right\Vert _{Y}^{2}+\left( 1-\frac{\lambda -\mu }{\lambda }\right) \left(
\left\Vert \nabla u\left( t\right) \right\Vert ^{2}+\left\Vert u\left(
t\right) \right\Vert _{L^{2}\left( \Gamma \right) }^{2}\right) \leq C,
\label{ineq-par}
\end{equation}%
for all $t\geq 0$, where we recall that $\mu \in \left( 0,\lambda \right] .$
Exploiting now the continuous embedding $\mathcal{V}^{1}\hookrightarrow Y$, (%
\ref{dissip-par}) follows from the application of Gronwall's inequality (see
Proposition \ref{GL}, Appendix) to (\ref{ineq-par}). The claim is proven.
\end{proof}

\begin{remark}
Theorem (\ref{t:parabolic-weak-solutions}) still holds if we keep (\ref%
{eq:assumption-f-3}) and we drop the assumptions (\ref{eq:assumption-f-1}), (%
\ref{eq:assumption-f-2}), and replace them by the following:%
\begin{equation}
\eta _{1}\left\vert y\right\vert ^{p}-C_{f}\leq f\left( y\right) y\leq \eta
_{2}\left\vert y\right\vert ^{p}+C_{f},  \label{n2}
\end{equation}%
for some $\eta _{1}$, $\eta _{2}>0,$ $C_{f}\geq 0$ and any $p>2$. In this
case, the same weak formulation (\ref{weakf-par})\ must be satisfied a.e. on
$\left[ 0,T\right] $, for all $\xi =(\chi ,\psi )\in \mathcal{V}^{1},$ with $%
\chi \in L^{p}\left( \Omega \right) $ (see, e.g., \cite{Gal12, Gal12-2}).
Finally, we note that without assumption (\ref{eq:assumption-f-3}), the
uniqueness of weak solutions (given in Definition \ref{weak-par}) is not
known in general (see \cite{Gal12-2}).
\end{remark}

\begin{corollary}
Let the assumptions of Theorem \ref{t:parabolic-weak-solutions} be
satisfied. We can define a strongly continuous semigroup%
\begin{equation*}
S_{0}\left( t\right) :Y\rightarrow Y
\end{equation*}%
by setting, for all $t\geq 0,$
\begin{equation*}
S_{0}(t)\zeta _{0}:=\zeta \left( t\right)
\end{equation*}%
where $\zeta \left( t\right) =(u\left( t\right) ,{u}_{\mid \Gamma }\left(
t\right) )$ is the unique weak solution to problem (\ref{eq:pde-p-1})-(\ref%
{eq:pde-p-3}).
\end{corollary}

The existence of a bounded absorbing set in $\mathcal{V}^{1}$ was shown for
the first time in \cite[Theorem 2.8]{Gal&Warma10} and the existence of the
global attractor for (\ref{eq:pde-p-1})-(\ref{eq:pde-p-3}) can be found in
\cite{Gal12-2}, \cite{Gal12}. The following theorem concerns the existence
and regularity of the global attractor $\mathcal{A}_{0}$ admitted by the
semiflow $S_{0}$ and is taken from \cite[Theorem 2.3]{Gal12}. The proof
relies on a uniform estimate which states that problem (\ref{eq:pde-p-1})-(%
\ref{eq:pde-p-3}) possesses the $Y-\mathcal{V}^{2}$ smoothing property and
exploits (\ref{dissip-par}).

\begin{theorem}
\label{t:global-attractor-and-regularity-p} The semiflow $S_{0}$ possesses a
connected global attractor $\mathcal{A}_{0}$ in $Y$, which is a bounded
subset of $\mathcal{V}^{2}$. The global attractor $\mathcal{A}_{0}$ contains
only strong solutions. Finally, $S_{0}$ also admits an exponential attractor
$\mathcal{M}_{0}$ which is bounded $\mathcal{V}^{2}$ and compact in $Y.$
\end{theorem}

\begin{remark}
The boundedness of $\mathcal{A}_{0}$ in $\mathcal{V}^{2}$, shown in \cite[%
Theorem 2.3]{Gal12}, is essential for the proof of the continuity property
at $\varepsilon =0$ of the global attractors associated with problem (\ref%
{eq:pde-h-1})-(\ref{eq:pde-h-3}). The last assertion follows from results in
\cite[Theorem 4.2]{GM09}, where (\ref{eq:pde-p-1})-(\ref{eq:pde-p-3}) is a
special case of a phase-field system endowed with dynamic boundary
conditions.
\end{remark}

\section{The hyperbolic relaxation problem}

\label{s:hyperbolic}

In this section, we study the hyperbolic relaxation problem (\ref{eq:pde-h-1}%
)-(\ref{eq:pde-h-3}) with $\varepsilon \in (0,1]$. Our first goal is to
prove the existence of a global attractor for (\ref{eq:pde-h-1})-(\ref%
{eq:pde-h-3}). As indicated in \cite{Wu&Zheng06}, semigroup methods are
applied to obtain local mild solutions whereby a suitable estimate is used
to extend the solution to a global one. We will offer a detailed
presentation on the well-posedness of the hyperbolic relaxation problem in
this section for the reader's convenience. The solution operators define a
semiflow on the phase space and because of the continuous dependence
estimate on the solutions, the semiflow is locally Lipschitz continuous,
uniformly in $t$ on compact intervals. Further estimates are used to
establish the existence of an absorbing set for the semiflow. As discussed
above, we will follow the decomposition method in \cite{Pata&Zelik06} to
obtain the existence of the global attractor in $H^1(\Omega)\times
L^2(\Omega)$ for the corresponding semiflow $S_{\varepsilon }$, for each $%
\varepsilon \in (0,1]$. The (optimal) regularity result for the global
attractors $\mathcal{A}_{\varepsilon }$\ and a proof of their continuity
properties conclude the section.

\subsection{The functional framework}

\label{s:notation}

Here, we consider the functional setup associated with problem (\ref%
{eq:pde-h-1})-(\ref{eq:pde-h-3}). The finite energy phase space for the
hyperbolic relaxation problem is the space
\begin{equation*}
\mathcal{H}_{\varepsilon }=H^{1}(\Omega )\times L^{2}(\Omega ).
\end{equation*}%
The space $\mathcal{H}_{\varepsilon }$ is Hilbert when endowed with the $%
\varepsilon $-weighted norm whose square is given by, for $\varphi =(u,v)\in
\mathcal{H}_{\varepsilon }=H^{1}(\Omega )\times L^{2}(\Omega )$,
\begin{equation*}
\Vert \varphi \Vert _{\mathcal{H}_{\varepsilon }}^{2}:=\Vert u\Vert
_{1}^{2}+\varepsilon \Vert v\Vert ^{2}=\left( \Vert \nabla u\Vert ^{2}+\Vert
u\Vert _{L^{2}(\Gamma )}^{2}\right) +\varepsilon \Vert v\Vert ^{2}.
\end{equation*}

As introduced in \cite{Wu&Zheng06} (cf. also \cite{CEL02}), $\Delta _{%
\mathrm{R}}:L^{2}(\Omega )\rightarrow L^{2}(\Omega )$ is the Robin-Laplacian
operator with domain
\begin{equation*}
D(\Delta _{\mathrm{R}})=\{u\in H^{2}(\Omega ):\partial _{\mathbf{n}}u+u=0~%
\text{on}~\Gamma \}.
\end{equation*}%
Easy calculations show that the operator $\Delta _{\mathrm{R}}$ is
self-adjoint and positive. The Robin-Laplacian is extended to a continuous
operator $\Delta _{\mathrm{R}}:H^{1}(\Omega )\rightarrow \left( H^{1}(\Omega
)\right) ^{\ast }$, defined by, for all $v\in H^{1}(\Omega )$,
\begin{equation*}
(-\Delta _{\mathrm{R}}u,v)=\langle \nabla u,\nabla v\rangle +\langle
u,v\rangle _{L^{2}(\Gamma )}.
\end{equation*}%
Next, \cite{CEL02,Wu&Zheng06} also define the Robin map $R:H^{s}(\Gamma
)\rightarrow H^{s+(3/2)}(\Omega )$ by
\begin{equation*}
Rp=q~\text{if and only if}~\Delta q=0~\text{in}~\Omega ,~\text{and}~\partial
_{\mathbf{n}}q+q=p~\text{on}~\Gamma .
\end{equation*}%
The adjoint of the Robin map satisfies, for all $v\in H^{1}(\Omega )$,
\begin{equation*}
R^{\ast }\Delta _{\mathrm{R}}v=-v~\text{on}~\Gamma .
\end{equation*}

Define the closed subspace of $H^{2}(\Omega )\times H^{1}(\Omega )$,
\begin{equation*}
\mathcal{D}_{\varepsilon }:=\{(u,v)\in H^{2}(\Omega )\times H^{1}(\Omega
):\partial _{\mathbf{n}}u+u=-v~\text{on}~\Gamma \}.
\end{equation*}%
endowed with norm whose square is given by, for all $\varphi =\left(
u,v\right) \in \mathcal{D}_{\varepsilon }$,
\begin{equation*}
\Vert \varphi \Vert _{\mathcal{D}_{\varepsilon }}^{2}:=\Vert u\Vert
_{2}^{2}+\Vert v\Vert _{1}^{2}.
\end{equation*}%
Let $D(A_{\varepsilon })=\mathcal{D}_{\varepsilon }$ (note that $\varepsilon
$-dependance does not enter through the norm of $\mathcal{D}_{\varepsilon },$
but rather in the definition of $A_{\varepsilon }$ below). Define the linear
unbounded operator $A_{\varepsilon }:D(A_{\varepsilon })\rightarrow \mathcal{%
H}_{\varepsilon }$ by
\begin{equation*}
A_{\varepsilon }:=%
\begin{pmatrix}
0 & 1 \\
\frac{1}{\varepsilon }\Delta _{\mathrm{R}} & \frac{1}{\varepsilon }(\Delta _{%
\mathrm{R}}R~\mathrm{tr_{D}}-1)%
\end{pmatrix}%
,
\end{equation*}%
where $\mathrm{tr_{D}}$ denotes the Dirichlet trace operator (i.e., $\mathrm{%
tr_{D}}(v)=v|_{\Gamma }$). Notice that if $(u,v)\in \mathcal{D}_{\varepsilon
}$, then $u+R\mathrm{tr_{D}}(v)\in D(\Delta _{\mathrm{R}})$. By the
Lumer-Phillips theorem (cf., e.g., \cite[Theorem I.4.3]{Pazy83}) and the
Lax-Milgram theorem, it is not hard to see that, for all $\varepsilon \in
(0,1]$, the operator $A_{\varepsilon }$, with domain $\mathcal{D}%
_{\varepsilon }$, is an infinitesimal generator of a strongly continuous
semigroup of contractions on $\mathcal{H}_{\varepsilon }$, denoted $%
e^{A_{\varepsilon }t}$.

Define the map $\mathcal{F}:\mathcal{H}_{\varepsilon }\rightarrow \mathcal{H}%
_{\varepsilon }$ by
\begin{equation*}
\mathcal{F}(\varphi ):=%
\begin{pmatrix}
0 \\
-\frac{1}{\varepsilon }f(u)%
\end{pmatrix}%
\end{equation*}%
for all $\varphi =(u,v)\in \mathcal{H}_{\varepsilon }$. Since $%
f:H^{1}(\Omega )\rightarrow L^{2}(\Omega )$ is locally Lipschitz continuous
\cite[cf., e.g., Theorem 2.7.13]{Zheng04}, it follows that the map $\mathcal{%
F}:\mathcal{H}_{\varepsilon }\rightarrow \mathcal{H}_{\varepsilon }$ is as
well.

The hyperbolic relaxation problem (\ref{eq:pde-h-1})-(\ref{eq:pde-h-3})\ may
be put into the abstract form in $\mathcal{H}_{\varepsilon }$, for $\varphi
(t)=(u(t),u_{t}(t))^{\mathrm{tr}}$,%
\begin{equation}
\displaystyle\frac{\mathrm{d}}{\mathrm{d}t}\varphi (t)=A_{\varepsilon
}\varphi (t)+\mathcal{F}(\varphi (t));~~\varphi (0)=%
\begin{pmatrix}
u_{0} \\
u_{1}%
\end{pmatrix}%
.  \label{eq:abstract-hyperbolic-problem}
\end{equation}

\begin{lemma}
For each $\varepsilon \in (0,1]$, the adjoint of $A_{\varepsilon }$, denoted
$A_{\varepsilon }^{\ast }$, is given by
\begin{equation*}
A_{\varepsilon }^{\ast }:=-%
\begin{pmatrix}
0 & 1 \\
\frac{1}{\varepsilon }\Delta _{\mathrm{R}} & -\frac{1}{\varepsilon }(\Delta
_{\mathrm{R}}R~\mathrm{tr_{D}}-1)%
\end{pmatrix}%
,
\end{equation*}%
with domain
\begin{equation*}
D(A_{\varepsilon }^{\ast }):=\{(\chi ,\psi )\in H^{2}(\Omega )\times
H^{1}(\Omega ):\partial _{\mathbf{n}}\chi +\chi =-\psi ~\text{on}~\Gamma \}.
\end{equation*}
\end{lemma}

\begin{proof}
The proof is a calculation similar to, e.g., \cite[Lemma 3.1]{Ball04}.
\end{proof}

\subsection{Well-posedness for the hyperbolic relaxation problem}

\label{well-posed}

The notion of weak solution to problem (\ref{eq:pde-h-1})-(\ref{eq:pde-h-3})
is as follows (see, \cite{Ball77}).

\begin{definition}
\label{mild}A function $\varphi =(u,u_{t}):[0,T]\rightarrow \mathcal{H}%
_{\varepsilon }$ is a weak solution of (\ref{eq:abstract-hyperbolic-problem}%
) on $[0,T],$ if and only if $\mathcal{F}(\varphi (\cdot ))\in L^{1}(0,T;%
\mathcal{H}_{\varepsilon })$ and $\varphi $ satisfies the variation of
constants formula, for all $t\in \lbrack 0,T]$,
\begin{equation*}
\varphi (t)=e^{A_{\varepsilon }t}\varphi _{0}+\int_{0}^{t}e^{A_{\varepsilon
}(t-s)}\mathcal{F}(\varphi (s))\mathrm{d}s.
\end{equation*}
\end{definition}

It can be easily shown that the notion of weak solution given in Definition %
\ref{mild} is also equivalent to the following notion of a weak solution
(see, e.g., \cite[Definition 3.1 and Proposition 3.5]{Ball04}).

\begin{definition}
\label{explicit}Let $T>0$ and $(u_{0},u_{1})\in \mathcal{H}_{\varepsilon }$.
A map $\varphi =(u,u_{t})\in C([0,T];\mathcal{H}_{\varepsilon })$ is a \emph{%
weak solution} of (\ref{eq:abstract-hyperbolic-problem}) on $[0,T],$ if for
each $\theta =(\chi ,\psi )\in D(A_{\varepsilon }^{\ast })$ the map $%
t\mapsto \langle \varphi (t),\theta \rangle _{\mathcal{H}_{\varepsilon }}$
is absolutely continuous on $[0,T]$ and satisfies, for almost all $t\in
\lbrack 0,T]$,
\begin{equation}
\frac{\mathrm{d}}{\mathrm{d}t}\langle \varphi (t),\theta \rangle _{\mathcal{H%
}_{\varepsilon }}=\langle \varphi (t),A_{\varepsilon }^{\ast }\theta \rangle
_{\mathcal{H}_{\varepsilon }}+\langle \mathcal{F}(\varphi (t)),\theta
\rangle _{\mathcal{H}_{\varepsilon }}.  \label{expf}
\end{equation}%
The map $\varphi =(u,u_{t})$ is a weak solution on $[0,\infty )$ (i.e., a
\emph{global weak solution}) if it is a weak solution on $[0,T]$, for all $%
T>0$.
\end{definition}

The above definitions are equivalent to the to the standard concept of a
weak (distributional) solution to (\ref{eq:pde-h-1})-(\ref{eq:pde-h-3}).

\begin{definition}
\label{weak}Let $\varepsilon \in (0,1]$. A function $\varphi
=(u,u_{t}):[0,T]\rightarrow \mathcal{H}_{\varepsilon }$ is a weak solution
of (\ref{eq:abstract-hyperbolic-problem}) (and, thus of (\ref{eq:pde-h-1})-(%
\ref{eq:pde-h-3})) on $[0,T],$ if%
\begin{equation*}
\varphi =(u,u_{t})\in C(\left[ 0,T\right] ;\mathcal{H}_{\varepsilon }),\text{
}u_{t}\in L^{2}(\left[ 0,T\right] \times \Gamma ),
\end{equation*}%
and, for each $\psi \in H^{1}\left( \Omega \right) ,$ $\left( u_{t},\psi
\right) \in C^{1}\left( \left[ 0,T\right] \right) $ with%
\begin{equation}
\frac{\mathrm{d}}{\mathrm{d}t}\left( \varepsilon u_{t}\left( t\right) ,\psi
\right) +\left\langle \nabla u\left( t\right) ,\nabla \psi \right\rangle
+\left\langle u_{t}\left( t\right) ,\psi \right\rangle +\left\langle
u_{t}\left( t\right) +u\left( t\right) ,\psi \right\rangle _{L^{2}\left(
\Gamma \right) }=-\left\langle f\left( u\left( t\right) \right) ,\psi
\right\rangle ,  \label{weakf}
\end{equation}%
for almost all $t\in \left[ 0,T\right] .$
\end{definition}

Indeed, by \cite[Lemma 3.3]{Ball04} we have that $f:H^{1}\left( \Omega
\right) \rightarrow L^{2}\left( \Omega \right) $ is sequentially weakly
continuous and continuous, on account of the assumptions (\ref%
{eq:assumption-f-1})-(\ref{eq:assumption-f-2}). Moreover, $\left( \varphi
_{t},\theta \right) \in C^{1}\left( \left[ 0,T\right] \right) $ for all $%
\theta \in D\left( A_{\varepsilon }^{\ast }\right) $, and (\ref{expf}) is
satisfied. The assertion in Definition \ref{weak} follows then from the
explicit characterization of $D\left( A_{\varepsilon }^{\ast }\right) $ and
from \cite[Proposition 3.4]{Ball04}.

Finally, the notion of strong solution to problem (\ref{eq:pde-h-1})-(\ref%
{eq:pde-h-3}) is as follows.

\begin{definition}
\label{d:regularity-property} \label{strong}Let $\varphi _{0}=\left(
u_{0},u_{1}\right) \in \mathcal{D}_{\varepsilon }$, $\varepsilon >0$, i.e., $%
(u_{0},u_{1})\in H^{2}(\Omega )\times H^{1}(\Omega )$ such that it satisfies
the compatibility condition%
\begin{equation*}
\partial _{\mathbf{n}}u_{0}+u_{0}+u_{1}=0,~\text{on}~\Gamma .
\end{equation*}%
A\textit{\ function }$\varphi \left( t\right) =\left( u\left( t\right)
,u_{t}\left( t\right) \right) $\textit{\ is called a (global) strong
solution if it is a weak solution in the sense of Definition \ref{weak}, and
if it satisfies the following regularity properties:}%
\begin{equation}
\begin{array}{l}
\varphi \in L^{\infty }(0,\infty ;\mathcal{D_{\varepsilon }})\text{, }%
\varphi_{t}\in L^{\infty }(0,\infty ;\mathcal{H}_{\varepsilon }), \\
~u_{tt}\in L^{\infty }(0,\infty ;L^{2}(\Omega )),\text{ }u_{tt}\in
L^{2}(0,\infty ;L^{2}(\Gamma )).%
\end{array}
\label{eq:regularity-property}
\end{equation}%
Therefore, $\varphi \left( t\right) =\left( u\left( t\right) ,u_{t}\left(
t\right) \right) $ satisfies the equations (\ref{eq:pde-h-1})-(\ref%
{eq:pde-h-3}) almost everywhere, i.e., is a strong solution.
\end{definition}

We can now state the main theorems of this section.

\begin{theorem}
\label{t:hyperbolic-weak-solutions} Assume (\ref{eq:assumption-f-1}) and (%
\ref{eq:assumption-f-2}) hold. For each $\varepsilon \in (0,1]$ and $\varphi
_{0}=(u_{0},u_{1})\in \mathcal{H}_{\varepsilon }$, there exists a unique
global weak solution $\varphi =(u,u_{t})\in C([0,\infty );\mathcal{H}%
_{\varepsilon })$ to (\ref{eq:pde-h-1})-(\ref{eq:pde-h-3}). In addition,
\begin{equation}
\partial _{\mathbf{n}}u\in L_{\mathrm{loc}}^{2}([0,\infty )\times \Gamma )~~%
\text{and}~~u_{t}\in L_{\mathrm{loc}}^{2}([0,\infty )\times \Gamma ).
\label{eq:bounded-boundary-1}
\end{equation}%
For each weak solution, the map
\begin{equation}
t\mapsto \Vert \varphi (t)\Vert _{\mathcal{H}_{\varepsilon
}}^{2}+2\int_{\Omega }F(u(t))\mathrm{d}x  \label{eq:C1-map}
\end{equation}%
is $C^{1}([0,\infty ))$ and the energy equation%
\begin{equation}
\frac{\mathrm{d}}{\mathrm{d}t}\left\{ \Vert \varphi (t)\Vert _{\mathcal{H}%
_{\varepsilon }}^{2}+2\int_{\Omega }F(u(t))\mathrm{d}x\right\} =-2\Vert
u_{t}(t)\Vert ^{2}-2\Vert u_{t}(t)\Vert _{L^{2}(\Gamma )}^{2}
\label{eq:energy-1}
\end{equation}%
holds (in the sense of distributions) a.e. on $[0,\infty )$. Furthermore,
let $\varphi (t)=(u(t),u_{t}(t))$ and $\theta (t)=(v(t),v_{t}(t))$ denote
the corresponding weak solution with initial data $\varphi
_{0}=(u_{0},u_{1})\in \mathcal{H}_{\varepsilon }$ and $\theta
_{0}=(v_{0},v_{1})\in \mathcal{H}_{\varepsilon }$, respectively, such that $%
\left\Vert \varphi _{0}\right\Vert _{\mathcal{H}_{\varepsilon }}\leq R,$ $%
\left\Vert \theta _{0}\right\Vert _{\mathcal{H}_{\varepsilon }}\leq R.$ Then
there exists a constant $\nu _{1}=\nu _{1}(R)>0$, such that, for all $t\geq
0 $,%
\begin{align}
& \Vert \varphi (t)-\theta (t)\Vert _{\mathcal{H}_{\varepsilon
}}^{2}+\int_{0}^{t}\left( \Vert u_{t}\left( \tau \right) -v_{t}\left( \tau
\right) \Vert ^{2}+\Vert u_{t}\left( \tau \right) -v_{t}\left( \tau \right)
\Vert _{L^{2}(\Gamma )}^{2}\right) \mathrm{d}\tau
\label{eq:continuous-dependence} \\
& \leq e^{\nu _{1}t}\Vert \varphi _{0}-\theta _{0}\Vert _{\mathcal{H}%
_{\varepsilon }}^{2}.  \notag
\end{align}
\end{theorem}

\begin{theorem}
\label{strong-sol-hyper}For each $\varepsilon \in (0,1]$ and $%
(u_{0},u_{1})\in \mathcal{D}_{\varepsilon }$, problem (\ref{eq:pde-h-1})-(%
\ref{eq:pde-h-3}) possesses a unique global strong solution in the sense of
Definition \ref{strong}.
\end{theorem}

\begin{remark}
The proof of Theorem \ref{strong-sol-hyper} is outlined in \cite{Wu&Zheng06}
(cf., also \cite{CEL02}) when $\varepsilon =1$.
\end{remark}

\begin{proof}[Proof of Theorem \protect\ref{t:hyperbolic-weak-solutions}]
We only give a sketch of the proof.

\emph{Step 1.} As discussed in the previous section, for each $\varepsilon
\in (0,1]$, the operator $A_{\varepsilon }$ with domain $D(A_{\varepsilon })=%
\mathcal{D}_{\varepsilon }$ is an infinitesimal generator of a strongly
continuous semigroup of contractions on $\mathcal{H}_{\varepsilon }$, and
the map $\mathcal{F}:\mathcal{H}_{\varepsilon }\rightarrow \mathcal{H}%
_{\varepsilon }$ is locally Lipschitz continuous. Therefore, by \cite[%
Theorem 2.5.4]{Zheng04}, for any $\varepsilon \in (0,1]$ and for any $%
\varphi _{0}=(u_{0},u_{1})\in \mathcal{H}_{\varepsilon }$, there is a $%
T^{\ast }=T^{\ast }(\Vert \varphi _{0}\Vert _{\mathcal{H}_{\varepsilon }})>0$%
, such that the abstract problem (\ref{eq:abstract-hyperbolic-problem})
admits a unique local weak solution on $[0,T^{\ast })$ satisfying
\begin{equation*}
\varphi \in C([0,T^{\ast });\mathcal{H}_{\varepsilon }).
\end{equation*}

The next step is to show that $T^{\ast }(\Vert \varphi _{0}\Vert _{\mathcal{H%
}_{\varepsilon }})=\infty $. Since the map (\ref{eq:C1-map}) is absolutely
continuous on $[0,T^{\ast })$ (cf., e.g., \cite[Theorem 3.1]{Ball04}), then
integration of the energy equation (\ref{eq:energy-1}) over $(0,t)$ yields,
for all $t\in \lbrack 0,T^{\ast })$,
\begin{equation}
\begin{aligned}\label{eq:hyperbolic-energy-4}
\|\varphi(t)\|^2_{\mathcal{H}_\varepsilon} & + 2\int_\Omega F(u(t)) {\rm{d}}
x + 2\int_0^t \|u_t(\tau)\|^2 {\rm{d}} \tau + 2\int_0^t
\|u_t(\tau)\|^2_{L^2(\Gamma)} {\rm{d}} \tau \\ & =
\|\varphi_0\|^2_{\mathcal{H}_\varepsilon} + 2\int_\Omega F(u_0) {\rm{d}} x.
\end{aligned}
\end{equation}%
Applying inequality (\ref{eq:consequence-F-1}) to (\ref%
{eq:hyperbolic-energy-4}) and applying (\ref{eq:consequence-F-2}) to the
integral on the right hand side, we find that there is a function $Q(\Vert
\varphi _{0}\Vert _{\mathcal{H}_{\varepsilon }})>0$, such that, for all $%
t\in \lbrack 0,T^{\ast })$,
\begin{equation}
\Vert \varphi (t)\Vert _{\mathcal{H}_\varepsilon }\leq Q(\Vert \varphi
_{0}\Vert _{\mathcal{H}_{\varepsilon }}).  \label{eq:hyperbolic-bound-1}
\end{equation}%
Since the bound on the right hand side of (\ref{eq:hyperbolic-bound-1}) is
independent of $t\in \lbrack 0,T^{\ast })$, $T^{\ast }(\Vert \varphi
_{0}\Vert _{\mathcal{H}_{\varepsilon }})$ can be extended indefinitely, and
therefore, for each $\varepsilon \in (0,1]$, we have that $T^{\ast }(\Vert
\varphi _{0}\Vert _{\mathcal{H}_{\varepsilon }})=\infty $.

We now show the boundary property (\ref{eq:bounded-boundary-1}). Applying (%
\ref{eq:consequence-F-1}), (\ref{eq:consequence-F-2}) and (\ref%
{eq:hyperbolic-bound-1}) to identity (\ref{eq:hyperbolic-energy-4}), we
obtain a bound of the form, for all $\varphi _{0}\in \mathcal{H}%
_{\varepsilon }$ and $t\geq 0$, in which
\begin{equation*}
\int_{0}^{t}\Vert u_{t}(\tau )\Vert _{L^{2}(\Gamma )}^{2}\mathrm{d}\tau \leq
Q(\Vert \varphi _{0}\Vert _{\mathcal{H}_{\varepsilon }}).
\end{equation*}%
It follows that $u_{t}\in L_{\mathrm{loc}}^{2}([0,\infty )\times \Gamma )$.
By the trace theorem, $u\in L^{\infty }(0,\infty ;H^{1}(\Omega
))\hookrightarrow L^{\infty }(0,\infty ;L^{2}(\Gamma ))$, so $u\in L_{%
\mathrm{loc}}^{2}([0,\infty )\times \Gamma )$. Comparison in (\ref%
{eq:pde-h-2}) yields that $\partial _{\mathbf{n}}u\in L_{\mathrm{loc}%
}^{2}([0,\infty )\times \Gamma )$.

\emph{Step 2.} To show that the continuous dependence estimate (\ref%
{eq:continuous-dependence}) holds, consider the difference $z(t):=u(t)-v(t)$%
, $t\geq 0$. We easily get%
\begin{equation}
\frac{\mathrm{d}}{\mathrm{d}t}\Vert (z,z_{t})\Vert _{\mathcal{H}%
_{\varepsilon }}^{2}+2\Vert z_{t}\Vert ^{2}+2\Vert z_{t}\Vert _{L^{2}(\Gamma
)}^{2}=2\langle f(v)-f(u),z_{t}\rangle .  \label{eq:dependence-2}
\end{equation}%
Since $f:H^{1}(\Omega )\rightarrow L^{2}(\Omega )$ is locally Lipschitz
continuous, then
\begin{equation}
2|\langle f(v)-f(u),z_{t}\rangle |\leq Q(R)\Vert z\Vert _{1}^{2}+\Vert
z_{t}\Vert ^{2},  \label{eq:dependence-3}
\end{equation}%
where $R>0$ is such that $\left\Vert \varphi _{0}\right\Vert _{\mathcal{H}%
_{\varepsilon }}\leq R,$ $\left\Vert \theta _{0}\right\Vert _{\mathcal{H}%
_{\varepsilon }}\leq R.$ Combining (\ref{eq:dependence-2}) and (\ref%
{eq:dependence-3}) produces, for almost all $t\geq 0$,
\begin{equation}
\frac{\mathrm{d}}{\mathrm{d}t}\Vert (z,z_{t})\Vert _{\mathcal{H}%
_{\varepsilon }}^{2}\leq Q(R)\Vert (z,z_{t})\Vert _{\mathcal{H}_{\varepsilon
}}^{2}.  \label{diff-ineq}
\end{equation}%
Hence, (\ref{eq:continuous-dependence}) follows immediately from (\ref%
{diff-ineq}), using the standard Gronwall lemma. This completes the proof of
the theorem.
\end{proof}

In view of Theorem \ref{t:hyperbolic-weak-solutions}, the following is
immediate.

\begin{corollary}
Let the assumptions of Theorem \ref{t:hyperbolic-weak-solutions} be
satisfied. Then, for each each $\varepsilon \in (0,1]$ we can define a
strongly continuous semigroup%
\begin{equation*}
S_{\varepsilon }\left( t\right) :\mathcal{H}_{\varepsilon }\rightarrow
\mathcal{H}_{\varepsilon },
\end{equation*}%
by setting, for all $t\geq 0,$
\begin{equation*}
S_{\varepsilon }\left( t\right) \varphi _{0}=\varphi \left( t\right) =\left(
u\left( t\right) ,u_{t}\left( t\right) \right) ,
\end{equation*}%
where $\varphi \left( t\right) $ is the unique weak solution to problem (\ref%
{eq:pde-h-1})-(\ref{eq:pde-h-3}).
\end{corollary}

\subsection{The global attractor $\mathcal{A}_{\protect\varepsilon }$ in $%
\mathcal{H}_{\protect\varepsilon }$}

\label{global-hyper}

In this section, we aim to show the existence of a global attractor, and
prove some additional regularity properties. We point out that all the
computations we will perform below can be rigorously justified by means of
an approximation procedure which relies upon the result in Theorem \ref%
{strong-sol-hyper}. Indeed, one shall use the usual procedure of
approximating weak solutions by strong solutions, and then pass to the limit
by using density theorems in the final estimates (see, also, \cite{CEL02}).
Thus, in what follows we can proceed formally.

We begin our analysis with a uniform estimate for the weak solutions of
Theorem \ref{t:hyperbolic-weak-solutions}. The estimate provides the
existence of a bounded absorbing set $\mathcal{B}_{\varepsilon }\subset
\mathcal{H}_{\varepsilon },$ for the semiflow $S_{\varepsilon }$, for each $%
\varepsilon \in (0,1]$.

\begin{lemma}
\label{t:bounded-absorbing-set-h} For all $\varepsilon \in (0,1]$ and $%
\varphi _{0}=(u_{0},u_{1})\in \mathcal{H}_{\varepsilon }$, there exist a
positive function $Q$, constants $\omega _{0}>0$, $P_{0}>0$, all independent
of $\varepsilon ,$ such that $\varphi (t)$ satisfies, for all $t\geq 0$,
\begin{equation}
\left\Vert \varphi \left( t\right) \right\Vert _{\mathcal{H}_{\varepsilon
}}^{2}\leq Q(\Vert \varphi _{0}\Vert _{\mathcal{H}_{\varepsilon
}})e^{-\omega _{0}t}+P_{0}.  \label{eq:decay-1}
\end{equation}%
Consequently, the ball $\mathcal{B}_{\varepsilon }$ in $\mathcal{H}%
_{\varepsilon }$,
\begin{equation}
\mathcal{B}_{\varepsilon }:=\{\varphi \in \mathcal{H}_{\varepsilon }:\Vert
\varphi \Vert _{\mathcal{H}_{\varepsilon }}\leq P_{0}+1\}  \label{absorbing}
\end{equation}%
is a bounded absorbing set in $\mathcal{H}_{\varepsilon }$ for the dynamical
system $\left( S_{\varepsilon }\left( t\right) ,\mathcal{H}_{\varepsilon
}\right) .$
\end{lemma}

\begin{proof}
Let $\varepsilon \in (0,1]$ and $\varphi _{0}=(u_{0},u_{1})\in \mathcal{H}%
_{\varepsilon }$. For $\alpha >0$ yet to be chosen, multiply (\ref%
{eq:pde-h-1}) by $\alpha u$ in $L^{2}(\Omega )$. Adding the result to the
energy equation (\ref{eq:energy-1}) above yields the differential identity,
which holds for almost all $t\geq 0$,
\begin{equation}
\begin{aligned} \frac{\diff}{\diff t} & \left\{
\|\varphi\|^2_{\mathcal{H}_\varepsilon} + \alpha\varepsilon\langle u_t,u
\rangle + 2\int_\Omega F(u) ~ {\rm{d}} x \right\} + \\ & +
(2-\varepsilon\alpha)\|u_t\|^2 + \alpha\langle u_t,u \rangle +
\alpha\|u\|^2_1 + \\ & + 2\|u_t\|^2_{L^2(\Gamma)} + \alpha\langle u_t,u
\rangle_{L^2(\Gamma)} + \alpha\langle f(u),u \rangle = 0.
\label{eq:energy-2} \end{aligned}
\end{equation}

For each $\varepsilon \in (0,1]$, define the functional, $E_{\varepsilon }:%
\mathcal{H}_{\varepsilon }\rightarrow \mathbb{R}$, by
\begin{equation}
E_{\varepsilon }(\varphi \left( t\right) )=\Vert \varphi \left( t\right)
\Vert _{\mathcal{H}_{\varepsilon }}^{2}+\alpha \varepsilon \langle
u_{t}\left( t\right) ,u\left( t\right) \rangle +2\int_{\Omega }F(u\left(
t\right) )\mathrm{d}x.  \label{eq:functional-E-1}
\end{equation}%
It is not hard to see that the map $t\mapsto E_{\varepsilon }(\varphi (t))$
is $C^{1}([0,\infty ))$; this essentially follows from equation (\ref%
{eq:C1-map}) of Theorem \ref{t:hyperbolic-weak-solutions}. First, we
estimate, for all $\eta >0$,
\begin{equation}
\alpha |\langle u_{t},u\rangle _{L^{2}(\Gamma )}|\leq \alpha \eta \Vert
u_{t}\Vert _{L^{2}(\Gamma )}^{2}+\frac{\alpha }{4\eta }\Vert u\Vert
_{L^{2}(\Gamma )}^{2},  \label{eq:estimate-1}
\end{equation}%
and with (\ref{eq:consequence-f-1}), we have,
\begin{equation}
\alpha |\langle f(u),u\rangle |\geq -\frac{\alpha (\lambda -\mu )}{\lambda }%
\Vert u\Vert _{1}^{2}-\alpha C.  \label{eq:estimate-2}
\end{equation}%
Combining (\ref{eq:energy-2}) with (\ref{eq:estimate-1})-(\ref{eq:estimate-2}%
) gives%
\begin{align}
& \frac{\mathrm{d}}{\mathrm{d}t}E_{\varepsilon }+(2-\alpha )\varepsilon
\Vert u_{t}\Vert ^{2}+\alpha \langle u_{t},u\rangle  \label{eq:estimate-5} \\
& +\alpha \left( 1-\frac{1}{4\eta }-\frac{\lambda -\mu }{\lambda }\right)
\Vert u\Vert _{1}^{2}+\left( 2-\alpha \eta \right) \Vert u_{t}\Vert
_{L^{2}(\Gamma )}^{2}  \notag \\
& \leq \alpha C  \notag
\end{align}%
Hence, for any $\eta >\frac{\lambda }{4\mu }$ and any $0<\alpha <\min \{2,%
\frac{2}{\eta }\}$, then $2-\eta >0$ and $2-\alpha \eta >0$,
\begin{equation*}
\omega _{0}:=\min \left\{ 2-\alpha ,\alpha \left( \frac{\mu }{\lambda }-%
\frac{1}{4\eta }\right) \right\} >0,
\end{equation*}%
and estimate (\ref{eq:estimate-5}) becomes, for almost all $t\geq 0$,
\begin{equation}
\frac{\mathrm{d}}{\mathrm{d}t}E_{\varepsilon }+\omega _{0}E_{\varepsilon
}+(2-\alpha \eta )\Vert u_{t}\Vert _{L^{2}(\Gamma )}^{2}\leq C_{\alpha }.
\label{eq:estimate-6}
\end{equation}%
Applying Gronwall's inequality (see, e.g., \cite[Lemma 5]{Pata&Zelik06}; cf.
also Proposition \ref{GL}, Appendix) to (\ref{eq:estimate-6}) produces, for
all $t\geq 0$,%
\begin{equation}
E_{\varepsilon }\left( \varphi \left( t\right) \right) \leq E_{\varepsilon
}\left( \varphi \left( 0\right) \right) e^{-\omega _{0}t}+C.
\label{eq:estimate-7}
\end{equation}

We now apply (\ref{eq:consequence-F-1}) to (\ref{eq:functional-E-1}) to
attain the bound,
\begin{equation}
E_{\varepsilon }\left( \varphi \right) \geq \varepsilon \left( 1-\frac{%
\alpha }{2}\right) \Vert u_{t}\Vert ^{2}+\left( 1-\frac{\alpha }{2\lambda }-%
\frac{\lambda -\mu }{\lambda }\right) \Vert u\Vert _{1}^{2}-C.
\label{eq:estimate-3}
\end{equation}%
After updating the smallness condition on $\alpha $ to $0<\alpha <\min \{2,%
\frac{2}{\eta },2\mu \}$, we see that for
\begin{equation*}
\omega _{1}:=\min \left\{ 1-\frac{\alpha }{2},1-\frac{\alpha }{2\lambda }-%
\frac{\lambda -\mu }{\lambda }\right\} >0,
\end{equation*}%
then, for all $t\geq 0$,
\begin{equation}
E_{\varepsilon }(\varphi (t))\geq \omega _{1}\Vert \varphi (t)\Vert _{%
\mathcal{H}_{\varepsilon }}^{2}-C.  \label{eq:E-bound-lower}
\end{equation}%
On the other hand, by estimating in a similar fashion, using (\ref%
{eq:consequence-F-1}), there holds for all $t\geq 0,$%
\begin{equation}
E_{\varepsilon }(\varphi (t))\leq Q\left( \Vert \varphi (t)\Vert _{\mathcal{H%
}_{\varepsilon }}\right) .  \label{eq:E-bound-upper}
\end{equation}%
Thus, estimate (\ref{eq:decay-1}) follows now from (\ref{eq:E-bound-lower}),
(\ref{eq:E-bound-upper}) and (\ref{eq:estimate-7}). The assertion (\ref%
{absorbing}) is an immediate consequence of (\ref{eq:decay-1}). This
concludes the proof.
\end{proof}

\begin{remark}
\label{r:bound} The following bounds are an immediate consequence of
estimate (\ref{eq:decay-1}):%
\begin{equation}
\limsup_{t\rightarrow \infty }\Vert \varphi (t)\Vert _{\mathcal{H}%
_{\varepsilon }}^{2}\leq P_{0}  \label{eq:uniform-bound-u}
\end{equation}%
and%
\begin{equation}
\int_{0}^{\infty }\left( \left\Vert u_{t}\left( \tau \right) \right\Vert
^{2}+\left\Vert u_{t}\left( \tau \right) \right\Vert _{L^{2}\left( \Gamma
\right) }^{2}\right) \mathrm{d} \tau \leq Q\left( \left\Vert \varphi
_{0}\right\Vert _{\mathcal{H}_{\varepsilon }}\right) .
\label{eq:uniform-bound-dissipation-integral}
\end{equation}%
The last bound is found by integrating the energy equation (\ref{eq:energy-1}%
) with respect to $t$ over $(0,\infty )$ and estimating the result with (\ref%
{eq:assumption-f-1}), (\ref{eq:assumption-f-2}), (\ref{eq:consequence-F-1}),
(\ref{eq:consequence-F-2}) and (\ref{eq:uniform-bound-u}).
\end{remark}

\begin{remark}
\label{r:rate-1} Note that the last assumption (\ref{eq:assumption-f-3})
(which is, $f^{^{\prime }}\left( s\right) \geq -\theta ,$ for all $s\in
\mathbb{R}$) is nowhere needed in the proofs of Theorem \ref%
{t:hyperbolic-weak-solutions}, Theorem \ref{strong-sol-hyper} (cf. \cite[%
Theorem 1.1 and Lemma 2.2]{Wu&Zheng06}) and Lemma \ref%
{t:bounded-absorbing-set-h}. It will only become important later (see (\ref%
{eq:beta}))\ when we establish the optimal regularity of the global
attractor for the hyperbolic problem (\ref{eq:pde-h-1})-(\ref{eq:pde-h-3}).
\end{remark}

The semiflow $S_{\varepsilon }$ admits a bounded absorbing set $\mathcal{B}%
_{\varepsilon }$ in $\mathcal{H}_{\varepsilon }$. To obtain a global
attractor, it suffices to prove that the semiflow admits a decomposition
into the sum of two operators, $S_{\varepsilon }=Z_{\varepsilon
}+K_{\varepsilon }$, where $Z_{\varepsilon }=(Z_{\varepsilon }(t))_{t\geq 0}$
and $K_{\varepsilon }=(K_{\varepsilon }(t))_{t\geq 0}$ are not necessarily
semiflows, but, operators that are uniformly decaying to zero, and uniformly
compact for large $t$, respectively. To obtain the compactness property for
the operator $K_{\varepsilon }$, recall that, when fractional powers of the
Laplacian are well-defined, one usually multiplies the PDE by the solution
and a suitable fractional power of the Laplacian; i.e., $(-\Delta )^{s}u$
for some $s>0$, then estimates using a stronger norm while keeping in mind
the uniform bound on $u$ and the null initial conditions. However, in our
case, the dynamic boundary condition does not allow us to proceed with the
usual argument to obtain the relative compactness of $K_{\varepsilon }$.
This is because the Laplacian equipped with the dynamic boundary condition (%
\ref{eq:pde-h-2}) is not self-adjoint nor positive. In turn, we cannot apply
the standard spectral theory to define fractional powers of the Laplacian.
So to obtain the relative compactness of $K_{\varepsilon }$, we follow the
approach in \cite{Pata&Zelik06}. The main tool is to differentiate the
equations with respect to time, and obtain uniform estimates for the new
equations. Such strategies also proved useful when dealing with a damped
wave equation with acoustic boundary conditions \cite{Frigeri10}, or a wave
equation with a nonlinear dynamic boundary condition \cite{CEL02}, and
hyperbolic relaxation of a Cahn-Hilliard equation with dynamic boundary
conditions \cite{CGG11}, \cite{GM12}.

Following an approach similar to the one taken in the above references,
first define
\begin{equation}
\psi (s):=f(s)+\beta s  \label{eq:beta}
\end{equation}%
for some constant $\beta \geq \vartheta $ to be determined later (in this
case, $\psi ^{\prime }(s)\geq 0$ thanks to assumption (\ref%
{eq:assumption-f-3})). Set $\Psi (s):=\int_{0}^{s}\psi (\sigma )\mathrm{d}%
\sigma $. Let $\varphi _{0}=(u_{0},u_{1})\in \mathcal{H}_{\varepsilon }$.
Then rewrite the hyperbolic relaxation problem into the system of equations
in $v$ and $w$, where $v+w=u$,
\begin{equation}
\left\{
\begin{array}{ll}
\varepsilon v_{tt}+v_{t}-\Delta v+\psi (u)-\psi (w)=0 & \text{in}~~(0,\infty
)\times \Omega , \\
\partial _{\mathbf{n}}v+v+v_{t}=0 & \text{on}~~(0,\infty )\times \Gamma , \\
v(0)=u_{0},~~v_{t}(0)=u_{1}+f\left( 0\right) -\beta u_{0} & \text{in}%
~~\Omega ,%
\end{array}%
\right.  \label{eq:pde-v}
\end{equation}%
and
\begin{equation}
\left\{
\begin{array}{ll}
\varepsilon w_{tt}+w_{t}-\Delta w+\psi (w)=\beta u & \text{in}~~(0,\infty
)\times \Omega , \\
\partial _{\mathbf{n}}w+w+w_{t}=0 & \text{on}~~(0,\infty )\times \Gamma , \\
w(0)=0,~~w_{t}(0)=-f\left( 0\right) +\beta u_{0} & \text{in}~~\Omega .%
\end{array}%
\right.  \label{eq:pde-w}
\end{equation}%
In view of Lemmas \ref{t:uniform-bound-w} and \ref{t:uniform-decay} below,
we define the one-parameter family of maps, $K_{\varepsilon }(t):\mathcal{H}%
_{\varepsilon }\rightarrow \mathcal{H}_{\varepsilon }$, by
\begin{equation*}
K_{\varepsilon }(t)\varphi _{0}:=\left( w(t),w_{t}(t)\right) ,
\end{equation*}%
where $\left( w,w_{t}\right) $ is a solution of (\ref{eq:pde-w}). With such $%
w$, we may define a second function $\left( v,v_{t}\right) $ as the solution
of (\ref{eq:pde-v}). Through the dependence of $v$ on $w$ and $\varphi
_{0}=(u_{0},u_{1})$, the solution of (\ref{eq:pde-v}) defines a
one-parameter family of maps, $Z_{\varepsilon }(t):\mathcal{H}_{\varepsilon
}\rightarrow \mathcal{H}_{\varepsilon }$, defined by
\begin{equation*}
Z_{\varepsilon }(t)\varphi _{0}:=\left( v(t),v_{t}(t)\right) .
\end{equation*}%
Notice that if $v$ and $w$ are solutions to (\ref{eq:pde-v}) and (\ref%
{eq:pde-w}), respectively, then the function $u:=v+w$ is a solution to the
original hyperbolic relaxation problem (\ref{eq:pde-h-1})-(\ref{eq:pde-h-3}).

The first lemma shows that the operators $K_{\varepsilon }$ are bounded in $%
\mathcal{H}_{\varepsilon }$, uniformly with respect to $\varepsilon $. The
result essentially follows from the existence of a bounded absorbing set $%
\mathcal{B}_{\varepsilon }$ in $\mathcal{H}_{\varepsilon }$ for $%
S_{\varepsilon }$ (recall (\ref{eq:uniform-bound-u})).

\begin{lemma}
\label{t:uniform-bound-w} Assume (\ref{eq:assumption-f-1}), (\ref%
{eq:assumption-f-2}) and (\ref{eq:assumption-f-3}) hold. For each $%
\varepsilon \in (0,1]$ and $\varphi _{0}=(u_{0},u_{1})\in \mathcal{H}%
_{\varepsilon }$, there exists a unique global weak solution $(w,w_{t})\in
C([0,\infty );\mathcal{H}_{\varepsilon })$ to problem (\ref{eq:pde-w})
satisfying
\begin{equation}
\partial _{\mathbf{n}}w\in L_{\mathrm{loc}}^{2}([0,\infty )\times \Gamma )~~%
\text{and}~~w_{t}\in L_{\mathrm{loc}}^{2}([0,\infty )\times \Gamma ).
\label{eq:bounded-boundary-2}
\end{equation}%
Moreover, for all $\varphi _{0}\in \mathcal{H}_{\varepsilon }$ with $%
\left\Vert \varphi _{0}\right\Vert _{\mathcal{H}_{\varepsilon }}\leq R$ for
all $\varepsilon \in (0,1]$, there holds for all $t\geq 0$,
\begin{equation}
\Vert K_{\varepsilon }(t)\varphi _{0}\Vert _{\mathcal{H}_{\varepsilon }}\leq
Q(R).  \label{eq:uniform-bound-w}
\end{equation}
\end{lemma}

The following result will be useful later on.

\begin{lemma}
\label{t:Gronwall-bound} For all $\varepsilon \in (0,1]$ and $\eta >0$,
there is a function $Q_{\eta }\left( \cdot \right) \sim \eta ^{-1}$, such
that for every $0\leq s\leq t$ and $\varphi _{0}=(u_{0},u_{1})\in \mathcal{B}%
_{\varepsilon }$,%
\begin{equation}
\int_{s}^{t}\left( \Vert w_{t}(\tau )\Vert ^{2}+\Vert u_{t}(\tau )\Vert
^{2}\right) \mathrm{d}\tau \leq \frac{\eta }{2}(t-s)+Q_{\eta }\left(
R\right) ,  \label{eq:Gronwall-bound-0}
\end{equation}%
where $R>0$ is such that $\left\Vert \varphi _{0}\right\Vert _{\mathcal{H}%
_{\varepsilon }}\leq R,$ for all $\varepsilon\in(0,1]$.
\end{lemma}

\begin{proof}
Let $\varepsilon \in (0,1]$ and $\varphi _{0}=(u_{0},u_{1})\in \mathcal{H}%
_{\varepsilon }$, with $\left\Vert \varphi _{0}\right\Vert _{\mathcal{H}%
_{\varepsilon }}\leq R.$ Adding the identity
\begin{equation*}
-2\beta \frac{\mathrm{d}}{\mathrm{d}t}\langle u,w\rangle =-2\beta \langle
u_{t},w\rangle -2\beta \langle u,w_{t}\rangle
\end{equation*}%
to equation%
\begin{equation}
\frac{\mathrm{d}}{\mathrm{d}t}\left\{ \Vert (w,w_{t})\Vert _{\mathcal{H}%
_{\varepsilon }}^{2}+2\int_{\Omega }\Psi (w)\mathrm{d}x\right\} +2\Vert
w_{t}\Vert ^{2}+\Vert w_{t}\Vert _{L^{2}(\Gamma )}^{2}=2\beta \langle
u,w_{t}\rangle  \label{eq:uniform-bound-w-2}
\end{equation}%
produces, for almost all $t\geq 0$,
\begin{equation}
\begin{aligned} \frac{\diff}{\diff t}\left\{
\|(w,w_t)\|^2_{\mathcal{H}_\varepsilon} + 2\int_\Omega \Psi(w) {\rm{d}} x -
2\beta\langle u,w \rangle \right\} & + 2\|w_t\|^2 + \|w_t\|^2_{L^2(\Gamma)}
\\ & = -2\beta\langle u_t,w \rangle. \end{aligned}
\label{eq:Gronwall-bound-1}
\end{equation}%
Using (\ref{eq:uniform-bound-w}), we estimate, for all $\eta >0$,
\begin{equation}
2\beta |\langle u_{t},w\rangle |\leq \eta +Q_{\eta }\left( R\right) \Vert
u_{t}\Vert ^{2}.  \label{eq:Gronwall-bound-3}
\end{equation}%
For each $\varepsilon \in (0,1]$, define the functional $W_{\varepsilon
}:H^{1}(\Omega )\times H^{1}(\Omega )\times L^{2}(\Omega )\rightarrow
\mathbb{R}$,
\begin{equation*}
W_{\varepsilon }(t):=\Vert (w(t),w_{t}(t))\Vert _{\mathcal{H}_{\varepsilon
}}^{2}+2\int_{\Omega }\Psi (w(t))\mathrm{d}x-2\beta \langle u(t),w(t)\rangle
.
\end{equation*}%
Because of (\ref{eq:consequence-F-2}), (\ref{eq:assumption-f-1}), (\ref%
{eq:assumption-f-2}), (\ref{eq:beta}), (\ref{eq:uniform-bound-u}) and (\ref%
{eq:uniform-bound-w}), we can easily check that for all $t\geq 0$ and $%
\varepsilon \in (0,1]$,
\begin{equation}
|W_{\varepsilon }(t)|\leq Q\left( R\right) .  \label{eq:Gronwall-bound-4}
\end{equation}%
We now combine (\ref{eq:Gronwall-bound-1}) and (\ref{eq:Gronwall-bound-3})
together as, for all $\eta >0$ and for almost all $t\geq 0$,
\begin{equation}
\frac{\mathrm{d}}{\mathrm{d}t}W_{\varepsilon }+2\Vert w_{t}\Vert ^{2}+\Vert
w_{t}\Vert _{L^{2}(\Gamma )}^{2}+2\Vert u_{t}\Vert ^{2}\leq \eta +\left(
Q_{\eta }\left( R\right) +2\right) \Vert u_{t}\Vert ^{2}.
\label{eq:Gronwall-bound-2}
\end{equation}%
Integrating (\ref{eq:Gronwall-bound-2}) over $(0,t)$, and recalling (\ref%
{eq:Gronwall-bound-4}), (\ref{eq:uniform-bound-dissipation-integral}), gives
the desired estimate in (\ref{eq:Gronwall-bound-0}). This proves the claim.
\end{proof}

The next result shows that the operators $Z_{\varepsilon }$ are uniformly
decaying to zero in $\mathcal{H}_{\varepsilon }$.

\begin{lemma}
\label{t:uniform-decay} Assume (\ref{eq:assumption-f-1}), (\ref%
{eq:assumption-f-2}) and (\ref{eq:assumption-f-3}) hold. For each $%
\varepsilon \in (0,1]$ and $\varphi _{0}=(u_{0},u_{1})\in \mathcal{H}%
_{\varepsilon }$, there exists a unique global weak solution $(v,v_{t})\in
C([0,\infty );\mathcal{H}_{\varepsilon })$ to problem (\ref{eq:pde-v})
satisfying
\begin{equation}
\partial _{\mathbf{n}}v\in L_{\mathrm{loc}}^{2}([0,\infty )\times \Gamma )~~%
\text{and}~~v_{t}\in L_{\mathrm{loc}}^{2}([0,\infty )\times \Gamma ).
\label{eq:bounded-boundary-3}
\end{equation}%
Moreover, for all $\varphi_0\in\mathcal{D}_\varepsilon$ with $\left\Vert
\varphi _{0}\right\Vert _{\mathcal{H}_{\varepsilon }}\leq R$ for all $%
\varepsilon \in (0,1]$, there exists $\omega >0 $, independent of $%
\varepsilon $, such that, for all $t\geq 0$,
\begin{equation}
\Vert Z_{\varepsilon }(t)\varphi _{0}\Vert _{\mathcal{H}_{\varepsilon }}\leq
Q(R)e^{-\omega t}.  \label{eq:uniform-decay}
\end{equation}
\end{lemma}

\begin{proof}
In a similar fashion to the arguments in Section \ref{well-posed}, the
existence of a global weak solution as well as (\ref{eq:bounded-boundary-3})
can be found. Because of (\ref{eq:uniform-bound-u}) and (\ref%
{eq:uniform-bound-w}), we know that the functions $(u(t),u_{t}(t))$ and $%
(w(t),w_{t}(t))$ are uniformly bounded in $\mathcal{H}_{\varepsilon }$ with
respect to $t$ and $\varepsilon $. It remains to show that (\ref%
{eq:uniform-decay}) holds.

Let $\varepsilon \in (0,1]$ and $\varphi _{0}=(u_{0},u_{1})\in \mathcal{H}%
_{\varepsilon }$, with $R>0$ such that $\left\Vert \varphi _{0}\right\Vert _{%
\mathcal{H}_{\varepsilon }}\leq R$. Observe that,
\begin{equation}
\begin{aligned}\label{eq:uniform-decay-1} 2\langle \psi(u)-\psi(w),v_t
\rangle = \frac{\diff}{\diff t} & \left\{ 2\langle \psi(u)-\psi(w),v \rangle
- \langle \psi'(u)v,v \rangle \right\} - \\ & - 2\langle
(\psi'(u)-\psi'(w))w_t,v \rangle + \langle \psi''(u)u_t,v^2 \rangle.
\end{aligned}
\end{equation}%
Multiply the first equation of (\ref{eq:pde-v}) by $2v_{t}+\alpha v$ in $%
L^{2}(\Omega )$, for $\alpha >0$ to be chosen later. We find that, with (\ref%
{eq:uniform-decay-1}),
\begin{equation}
\begin{aligned} \frac{\diff}{\diff t} & \left\{ \varepsilon\|v_t\|^2 +
\alpha\varepsilon\langle v_t,v \rangle + \|v\|^2_1 + 2\langle
\psi(u)-\psi(w),v \rangle - \langle \psi'(u)v,v \rangle \right\} + \\ & +
(2-\alpha\varepsilon)\|v_t\|^2 + \alpha\langle v_t,v \rangle +
\alpha\|v\|^2_1 + 2\|v_t\|^2_{L^2(\Gamma)} + \alpha\langle v_t,v
\rangle_{L^2(\Gamma)} + \\ & + \alpha\langle \psi(u)-\psi(w),v \rangle =
2\langle (\psi'(u)-\psi'(w))w_t,v \rangle - \langle \psi''(u)u_t,v^2
\rangle. \label{eq:zero-decay-1} \end{aligned}
\end{equation}

For each $\varepsilon \in (0,1]$, define the functional%
\begin{equation*}
V_{\varepsilon }:H^{1}(\Omega )\times H^{1}(\Omega )\times H^{1}(\Omega
)\times L^{2}(\Omega )\rightarrow \mathbb{R},
\end{equation*}%
by%
\begin{equation*}
\begin{aligned} V_\varepsilon(t):= & \varepsilon\|v_t(t)\|^2 +
\alpha\varepsilon\langle v_t(t),v(t) \rangle + \|v(t)\|^2_1 + \\ & +
2\langle \psi(u(t))-\psi(w(t)),v(t) \rangle - \langle \psi'(u(t))v(t),v(t)
\rangle. \end{aligned}
\end{equation*}%
As with the functional $E_{\varepsilon }$ above, the map $t\mapsto
V_{\varepsilon }(t)$ is $AC(\mathbb{R}_{\geq 0};\mathbb{R}_{\geq 0})$. We
now will show that, given $(u,u_{t}),(w,w_{t})\in \mathcal{H}_{\varepsilon }$
are uniformly bounded with respect to $t$ and $\varepsilon $, there are
constants, $C_{1},C_{2}>0$, independent of $t$ and $\varepsilon $ (possibly
depending on $R>0$), in which, for all $(v,v_{t})\in \mathcal{H}%
_{\varepsilon }$,%
\begin{equation}
C_{1}\Vert (v,v_{t})\Vert _{\mathcal{H}_{\varepsilon }}^{2}\leq
V_{\varepsilon }\leq C_{2}\Vert (v,v_{t})\Vert _{\mathcal{H}_{\varepsilon
}}^{2}.  \label{assy}
\end{equation}%
We begin by estimating the products in $V_{\varepsilon }$ that involve $\psi
$; with (\ref{eq:assumption-f-1}), (\ref{eq:assumption-f-2}), the embedding $%
H^{1}(\Omega )\hookrightarrow L^{6}(\Omega )$ and (\ref{eq:uniform-bound-u}%
), there holds%
\begin{equation}
\begin{aligned} |\langle \psi ^{\prime }(u)v,v\rangle | & \leq C_{\Omega
}\left( 1+\Vert u\Vert _{1}^{2}\right) \Vert v\Vert _{1}\Vert v\Vert \\ &
\leq \frac{1}{2}\Vert v\Vert _{1}^{2}+Q(R)\Vert v\Vert ^{2}.
\label{eq:zero-decay-5} \end{aligned}
\end{equation}%
From assumption (\ref{eq:assumption-f-3}) and the definition of $\psi ,$ cf.
(\ref{eq:beta}),
\begin{equation}
2\langle \psi (u)-\psi (w),v\rangle \geq 2(\beta -\vartheta )\Vert v\Vert
^{2}.  \label{eq:zero-decay-51}
\end{equation}%
Hence, for $\beta $ sufficiently large, $\beta \geq \left( C\left( R\right)
+2\vartheta \right) /2$, the combination of (\ref{eq:zero-decay-5}) and (\ref%
{eq:zero-decay-51}) produces,
\begin{equation*}
\begin{aligned} 2\langle \psi(u)-\psi(w),v \rangle - \langle \psi'(u)v,v
\rangle & \geq 2(\beta-\vartheta)\|v\|^2-\frac{1}{2}\|v\|^2_1-C(R)\|v\|^2 \\
& \geq -\frac{1}{2}\|v\|^2_1. \end{aligned}
\end{equation*}%
Then we attain the lower bound on $V_{\varepsilon }$,
\begin{equation*}
V_{\varepsilon }\geq \left( 1-\frac{\alpha }{2}\right) \varepsilon \Vert
v_{t}\Vert ^{2}+\left( \frac{1}{2}-\frac{\alpha }{2\lambda }\right) \Vert
v\Vert _{1}^{2}.
\end{equation*}%
So for $0<\alpha <\min \{2,\lambda \}$, set
\begin{equation*}
\omega _{2}:=\min \left\{ 1-\frac{\alpha }{2},\frac{1}{2}-\frac{\alpha }{%
2\lambda }\right\} >0,
\end{equation*}%
then, for all $t\geq 0$, we have that
\begin{equation}
V_{\varepsilon }(t)\geq \omega _{2}\Vert (v(t),v_{t}(t))\Vert _{\mathcal{H}%
_{\varepsilon }}^{2}.  \label{eq:zero-decay-8}
\end{equation}%
Now by the (local) Lipschitz continuity of $f,$ and the uniform bounds on $u$
and $w$, it is easy to check that%
\begin{equation*}
2\langle \psi (u)-\psi (w),v\rangle \leq 2\Vert \psi (u)-\psi (w)\Vert \Vert
v\Vert \leq Q(R)\Vert v\Vert _{1}^{2}.
\end{equation*}%
Also, using (\ref{eq:assumption-f-1}), (\ref{eq:assumption-f-2}) and the
bound (\ref{eq:uniform-bound-u}), there also holds%
\begin{equation}
|\langle \psi ^{\prime }(u)v,v\rangle |\leq Q(R)\Vert v\Vert _{1}^{2}.
\label{eq:zero-decay-11}
\end{equation}%
Thus, the assertion in (\ref{assy}) holds. Exploiting the fact that%
\begin{equation*}
\alpha |\langle v_{t},v\rangle _{L^{2}(\Gamma )}|\leq \frac{\alpha }{2}\Vert
v_{t}\Vert _{L^{2}(\Gamma )}^{2}+\frac{\alpha }{2}\Vert v\Vert
_{L^{2}(\Gamma )}^{2},
\end{equation*}%
we see that (\ref{eq:zero-decay-1}) becomes%
\begin{equation}
\begin{aligned} \frac{\diff}{\diff t} & V_\varepsilon +
(2-\alpha)\varepsilon\|v_t\|^2 + \alpha\langle v_t,v \rangle +
\alpha\|\nabla v\|^2 + \frac{\alpha}{2}\|v\|^2_{L^2(\Gamma)} + \\ & +
(2-\frac{\alpha}{2})\|v_t\|^2_{L^2(\Gamma)} + \alpha\langle
\psi(u)-\psi(w),v \rangle - \langle \psi'(u)v,v \rangle \\ & \leq -\langle
\psi'(u)v,v \rangle + 2\langle (\psi'(u)-\psi'(w))w_t,v \rangle - \langle
\psi''(u)u_t,v^2 \rangle. \label{eq:zero-decay-3} \end{aligned}
\end{equation}%
Recall that $0<\alpha <\min \{2,\lambda \}$, so when we set
\begin{equation*}
\omega _{3}:=\min \left\{ 2-\alpha ,1,\frac{\alpha }{2}\right\} >0,
\end{equation*}%
we write (\ref{eq:zero-decay-3}) as
\begin{equation}
\begin{aligned} \frac{\diff}{\diff t} & V_\varepsilon + \omega_3
V_\varepsilon \leq -\langle \psi'(u)v,v \rangle + 2\langle
(\psi'(u)-\psi'(w))w_t,v \rangle - \langle \psi''(u)u_t,v^2 \rangle.
\label{eq:zero-decay-4} \end{aligned}
\end{equation}%
Using the uniform bound on $u$ and $w$ (recall assumptions (\ref%
{eq:assumption-f-1}), (\ref{eq:assumption-f-2}), (\ref{eq:uniform-bound-u})
and (\ref{eq:uniform-bound-w})), there is a positive function $Q_{\eta
}(R)>0 $, depending on $\eta $, such that, for all $\eta >0$,%
\begin{align}
\left\vert \left\langle \left( \psi ^{\prime }(u)-\psi ^{^{\prime }}\left(
w\right) \right) w_{t},v\right\rangle \right\vert & \leq C_{\Omega }\left(
1+\left\Vert u\right\Vert _{1}+\left\Vert w\right\Vert _{1}\right)
\left\Vert w_{t}\right\Vert \left\Vert v\right\Vert _{1}^{2}
\label{eq:zero-decay-6} \\
& \leq \frac{\eta }{2}\left\Vert v\right\Vert _{1}^{2}+Q_{\eta }\left(
R\right) \left\Vert w_{t}\right\Vert ^{2}V_{\varepsilon }.  \notag
\end{align}%
The last inequality in the above estimate follows from (\ref{eq:zero-decay-8}%
). In a similar fashion we estimate using assumption (\ref{eq:assumption-f-1}%
) and the bound (\ref{eq:uniform-bound-u}),%
\begin{align}
|\langle \psi ^{^{\prime \prime }}(u)u_{t},v^{2}\rangle |& \leq C_{\Omega
}\left( 1+\left\Vert u\right\Vert _{1}\right) \left\Vert u_{t}\right\Vert
\left\Vert v\right\Vert _{1}^{2}  \label{eq:zero-decay-7} \\
& \leq \frac{\eta }{2}\left\Vert v\right\Vert _{1}^{2}+Q_{\eta }\left(
R\right) \left\Vert u_{t}\right\Vert ^{2}V_{\varepsilon }.  \notag
\end{align}%
Applying (\ref{eq:zero-decay-5}) to (\ref{eq:zero-decay-4}) and inserting (%
\ref{eq:zero-decay-6}) and (\ref{eq:zero-decay-7}) into (\ref%
{eq:zero-decay-4}), we then have%
\begin{equation}
\frac{\mathrm{d}}{\mathrm{d}t}V_{\varepsilon }+\omega _{3}V_{\varepsilon
}-\eta \Vert v\Vert _{1}^{2}\leq Q_{\eta }\left( R\right) \left( \Vert
u_{t}\Vert ^{2}+\Vert w_{t}\Vert ^{2}\right) V_{\varepsilon }.
\label{eq:zero-decay-9}
\end{equation}%
There is a sufficiently small $\eta $, precisely, $0<\eta <\omega _{3}/2$,
so that (\ref{eq:zero-decay-9}) becomes
\begin{equation}
\frac{\mathrm{d}}{\mathrm{d}t}V_{\varepsilon }+\eta V_{\varepsilon }\leq
Q_{\eta }\left( R\right) \left( \Vert u_{t}\Vert ^{2}+\Vert w_{t}\Vert
^{2}\right) V_{\varepsilon }.  \label{eq:zero-decay-10}
\end{equation}%
At this point, we remind the reader of Lemma \ref{t:Gronwall-bound}.
Applying a suitable Gronwall type inequality (see, e.g., \cite[Lemma 5]%
{Pata&Zelik06}; cf. also Proposition \ref{GL}, Appendix) to (\ref%
{eq:zero-decay-10}) yields%
\begin{equation}
V_{\varepsilon }(t)\leq V_{\varepsilon }(0)e^{Q_{\eta }\left( R\right)
}e^{-\eta t/2}.  \label{eq:zero-decay-10bis}
\end{equation}%
By virtue of (\ref{assy}), for all $\varepsilon \in (0,1],$%
\begin{eqnarray*}
V_{\varepsilon }(0) &\leq &Q\left( R\right) \Vert (v\left( 0\right)
,v_{t}\left( 0\right) )\Vert _{\mathcal{H}_{\varepsilon }}^{2} \\
&\leq &Q\left( R\right) \left( \left\Vert u_{0}\right\Vert
_{1}^{2}+\varepsilon \left\Vert u_{1}+f\left( 0\right) -\beta
u_{0}\right\Vert ^{2}\right) \\
&\leq &Q\left( R\right) ,
\end{eqnarray*}%
for some positive function $Q$ independent of $\varepsilon $. Therefore (\ref%
{eq:zero-decay-10bis}) shows that the operators $Z_{\varepsilon }$ are
uniformly decaying to zero.
\end{proof}

The following lemma establishes the uniform compactness of the operators $%
K_{\varepsilon }$.

\begin{lemma}
\label{t:uniform-compactness} For all $\varphi _{0}=(u_{0},u_{1})\in
\mathcal{H}_{\varepsilon }$ such that $\left\Vert \varphi _{0}\right\Vert _{%
\mathcal{H}_{\varepsilon }}\leq R$ for all $\varepsilon \in (0,1]$, the
following estimate holds:
\begin{equation*}
\Vert K_{\varepsilon }(t)\varphi _{0}\Vert _{\mathcal{D}_{\varepsilon }}\leq
Q(R),
\end{equation*}%
for all $t\geq 0.$ Furthermore, the operators $K_{\varepsilon }$ are
uniformly compact in $\mathcal{H}_{\varepsilon }$.
\end{lemma}

\begin{proof}
Let $\varepsilon \in (0,1]$ and let $\varphi _{0}=(u_{0},u_{1})\in \mathcal{H%
}_{\varepsilon }$ with $R>0$ such that $\left\Vert \varphi _{0}\right\Vert _{%
\mathcal{H}_{\varepsilon }}\leq R$. Differentiate (\ref{eq:pde-w}) with
respect to $t$ and set $h=w_{t}.$ Then $h$ satisfies the equations%
\begin{equation}
\left\{
\begin{array}{ll}
\varepsilon h_{tt}+h_{t}-\Delta h+\psi ^{\prime }(w)h=\beta u_{t} & \text{in}%
~~(0,\infty )\times \Omega , \\
\partial _{\mathbf{n}}h+h+h_{t}=0, & \text{on}~~(0,\infty )\times \Gamma ,
\\
h(0)=w_{t}\left( 0\right) ,~~h_{t}(0)=w_{tt}\left( 0\right) & \text{in}%
~~\Omega .%
\end{array}%
\right.  \label{eq:pde-diff-w}
\end{equation}%
Note that, by the choice of data in (\ref{eq:pde-w}), we actually have $%
h\left( 0\right) =-f\left( 0\right) +\beta u_{0}$ and $h_{t}\left( 0\right)
=0.$ Multiply the first equation of (\ref{eq:pde-diff-w}) by $2h_{t}+\alpha
h $, where $\alpha >0$ is yet to be determined, and integrate over $\Omega $%
. Adding the result to the identity
\begin{equation*}
2\langle \psi ^{\prime }(w)h,h_{t}\rangle =\frac{\mathrm{d}}{\mathrm{d}t}%
\langle \psi ^{\prime }(w)h,h\rangle -\langle \psi ^{\prime \prime
}(w)w_{t},h^{2}\rangle
\end{equation*}%
produces%
\begin{equation}
\begin{aligned} \frac{\diff}{\diff t} & \left\{ \varepsilon\|h_t\|^2 +
\alpha\varepsilon\langle h_t,h \rangle + \|h\|^2_1 + \langle \psi'(w)h,h
\rangle \right\} + \\ & + (2-\alpha\varepsilon)\|h_t\|^2 + \alpha\langle
h_t,h \rangle + \alpha\|h\|^2_1 + \\ & + 2\|h_t\|^2_{L^2(\Gamma)} +
\alpha\langle h_t,h \rangle_{L^2(\Gamma)} + \alpha\langle \psi'(w)h,h
\rangle \\ & = \langle \psi''(w)w_t,h^2 \rangle + 2\beta\langle u_t,h_t
\rangle + \alpha\beta \langle u_t,w_t \rangle. \label{eq:diff-w-1}
\end{aligned}
\end{equation}%
For each $\varepsilon \in (0,1]$, define the functional%
\begin{equation*}
\Psi _{\varepsilon }:H^{1}(\Omega )\times H^{1}(\Omega )\times L^{2}(\Omega
)\rightarrow \mathbb{R}
\end{equation*}%
by%
\begin{equation}
\Psi _{\varepsilon }(t):=\varepsilon \Vert h_{t}(t)\Vert ^{2}+\alpha
\varepsilon \langle h_{t}(t),h(t)\rangle +\Vert h(t)\Vert _{1}^{2}+\langle
\psi ^{\prime }(w(t))h(t),h(t)\rangle .  \label{eq:diff-w-2}
\end{equation}%
The map $t\mapsto \Psi _{\varepsilon }(t)$ is $AC(\mathbb{R}_{\geq 0};%
\mathbb{R}_{\geq 0})$. Because of the bound given in Lemma \ref%
{t:uniform-bound-w}, we obtain the estimate similar to (\ref%
{eq:zero-decay-11}),
\begin{equation}
\alpha |\langle \psi ^{\prime }(w)h,h\rangle |\leq \alpha Q(R)\Vert h\Vert
_{1}^{2}.  \label{eq:diff-w-3}
\end{equation}%
Obviously, we have%
\begin{equation}
\alpha \varepsilon |\langle h_{t},h\rangle |\leq \frac{\alpha \varepsilon }{2%
}\Vert h_{t}\Vert ^{2}+\frac{\alpha }{2\lambda }\Vert h\Vert _{1}^{2}.
\label{eq:diff-w-10}
\end{equation}%
After combining (\ref{eq:diff-w-2})-(\ref{eq:diff-w-10}), we find
\begin{equation*}
\Psi _{\varepsilon }\geq \left( 1-\frac{\alpha }{2}\right) \varepsilon \Vert
h_{t}\Vert ^{2}+\left( 1-\frac{\alpha }{2\lambda }-\alpha Q(R)\right) \Vert
h\Vert _{1}^{2}.
\end{equation*}%
Hence, when
\begin{equation*}
0<\alpha <\min \left\{ 2,\left( \frac{1}{2\lambda }+Q(R)\right)
^{-1}\right\} ,
\end{equation*}%
then,
\begin{equation*}
\omega _{4}(R):=\min \left\{ 1-\frac{\alpha }{2},1-\frac{\alpha }{2\lambda }%
-\alpha Q(R)\right\} >0,
\end{equation*}%
thus, for all $t\geq 0$
\begin{equation}
\Psi _{\varepsilon }(t)\geq \omega _{4}\Vert (h(t),h_{t}(t))\Vert _{\mathcal{%
H}_{\varepsilon }}^{2}.  \label{eq:diff-lower}
\end{equation}%
On the other hand, again with (\ref{eq:diff-w-3}),
\begin{equation*}
\Psi _{\varepsilon }\leq \left( 1+\frac{\alpha }{2}\right) \varepsilon \Vert
h_{t}\Vert ^{2}+\left( 1+\frac{\alpha }{2\lambda }+\alpha Q(R)\right) \Vert
h\Vert _{1}^{2},
\end{equation*}%
and with
\begin{equation*}
\omega _{5}(R):=\max \left\{ 1+\frac{\alpha }{2},1+\frac{\alpha }{2\lambda }%
+\alpha Q(R)\right\} ,
\end{equation*}%
an upper-bound for $\Psi _{\varepsilon }$ is given by, for all $t\geq 0$,
\begin{equation}
\Psi _{\varepsilon }(t)\leq \omega _{5}\Vert (h(t),h_{t}(t))\Vert _{\mathcal{%
H}_{\varepsilon }}^{2}.  \label{eq:diff-upper}
\end{equation}%
Using the bounds found in (\ref{eq:uniform-bound-u}) and (\ref%
{eq:uniform-bound-w}), we estimate the following terms from (\ref%
{eq:diff-w-1}), for all $\eta >0$,
\begin{equation}
\alpha |\langle h_{t},h\rangle _{L^{2}(\Gamma )}|\leq \alpha \eta \Vert
h_{t}\Vert _{L^{2}(\Gamma )}^{2}+\frac{\alpha }{4\eta }\Vert h\Vert
_{L^{2}(\Gamma )}^{2},  \label{eq:diff-w-4}
\end{equation}%
and
\begin{equation}
\begin{aligned}\label{eq:diff-w-5} 2\beta|\langle u_t,h_t \rangle| +
\alpha\beta|\langle u_t,w_t \rangle| & \leq Q(R)\|h_t\|+Q(R) \\ & \leq
\eta\|h_t\|^2+Q_\eta(R). \end{aligned}
\end{equation}%
Also, similar to (\ref{eq:zero-decay-7}), but when we now employ (\ref%
{eq:diff-lower}), we have that, for all $\eta >0$,
\begin{equation}
\langle \psi ^{\prime \prime }(w)w_{t},h^{2}\rangle \leq Q_{\eta }(R)\Vert
w_{t}\Vert \Psi _{\varepsilon }.  \label{eq:diff-w-6}
\end{equation}%
Combine (\ref{eq:diff-w-4})-(\ref{eq:diff-w-6}) with (\ref{eq:diff-w-1}) and
obtain the following estimate (note that when $2-\alpha -\eta >0$, we have $%
(2-\alpha -\eta )\varepsilon <2-\alpha \varepsilon -\eta $):%
\begin{equation}
\begin{aligned} \frac{\diff}{\diff t} & \Psi_\varepsilon +
(2-\alpha-\eta)\varepsilon\|h_t\|^2 + \alpha\langle h_t,h \rangle +
\alpha\|\nabla h\|^2 + \alpha\left( 1-\frac{1}{4\eta}
\right)\|h\|^2_{L^2(\Gamma)} + \\ & + (2-\alpha\eta)\|h_t\|^2_{L^2(\Gamma)}
+ \alpha\langle \psi'(w)h,h \rangle \leq Q_\eta(R)\|w_t\|\Psi_\varepsilon +
Q_\eta(R). \label{eq:diff-w-7} \end{aligned}
\end{equation}%
With some $\frac{1}{4}<\eta <2$ now fixed, then, for
\begin{equation*}
0<\alpha <\min \left\{ 2-\eta ,1,\frac{2}{\eta }\right\} \text{ and }\omega
_{6}:=1-\frac{1}{4\eta },
\end{equation*}%
we have%
\begin{equation}
\frac{\mathrm{d}}{\mathrm{d}t}\Psi _{\varepsilon }+\omega _{6}\Psi
_{\varepsilon }+(2-\alpha \eta )\Vert h_{t}\Vert _{L^{2}(\Gamma )}^{2}\leq
Q\left( R\right) \Vert w_{t}\Vert \Psi _{\varepsilon }+Q(R).
\label{eq:diff-w-20}
\end{equation}%
An immediate consequence of (\ref{eq:uniform-bound-dissipation-integral}) is
the bound on the following integral%
\begin{equation*}
\int_{0}^{\infty }\Vert w_{t}(\tau )\Vert ^{2}\mathrm{d}\tau \leq Q\left(
R\right) .
\end{equation*}%
Applying a suitable version of the Gronwall inequality (see, e.g., \cite[%
Lemma 2.2]{GP}; cf. also Proposition \ref{GL}, Appendix) it follows that
\begin{equation}
\Psi _{\varepsilon }(t)\leq Q\left( R\right) \Psi _{\varepsilon
}(0)e^{-\omega _{6}t/2}+Q\left( R\right) .  \label{eq:for-C1}
\end{equation}%
Using (\ref{eq:diff-lower}) and (\ref{eq:diff-upper}), and the fact that $%
\Psi _{\varepsilon }(0)\leq \omega _{5}\Vert (h(0),h_{t}(0))\Vert _{\mathcal{%
H}_{\varepsilon }}^{2}\leq Q\left( R\right) $, we arrive at the bound%
\begin{equation}
\Vert w_{t}(t)\Vert _{1}^{2}+\varepsilon \Vert w_{tt}(t)\Vert ^{2}\leq
Q\left( R\right) ,  \label{eq:wt-bounds}
\end{equation}%
for all $t\geq 0$, $\varepsilon \in (0,1]$ and $\varphi _{0}\in \mathcal{H}%
_{\varepsilon }$, with $R>0$ such that $\left\Vert \varphi _{0}\right\Vert _{%
\mathcal{H}_{\varepsilon }}\leq R.$

In order to bound $\Vert (w,w_{t})\Vert _{\mathcal{D}_{\varepsilon }}$, we
need to bound the term $\Vert w\Vert _{2}$. We have owing to standard
elliptic regularity theory (see, e.g., \cite[Theorem II.5.1]{Lions&Magenes72}%
), that%
\begin{equation}
\Vert w(t)\Vert _{2}\leq C\left( \Vert \Delta w(t)\Vert +\Vert \partial _{%
\mathbf{n}}w(t)\Vert _{H^{1/2}(\Gamma )}\right) .  \label{eq:H2-regularity}
\end{equation}%
Thus, using the first equation of (\ref{eq:pde-w})$,$ the bounds (\ref%
{eq:uniform-bound-u}), (\ref{eq:uniform-bound-w}) and (\ref{eq:wt-bounds}),
and also (\ref{eq:assumption-f-1}), (\ref{eq:assumption-f-2}) and (\ref%
{eq:beta}), we have
\begin{equation}
\Vert \Delta w\left( t\right) \Vert \leq \sqrt{\varepsilon }\Vert
w_{tt}\left( t\right) \Vert +\Vert w_{t}\left( t\right) \Vert +\Vert \psi
(w\left( t\right) )\Vert +\beta \Vert u\left( t\right) \Vert \leq Q(R).
\label{eq:wt-bounds-1}
\end{equation}%
Also, by (\ref{eq:wt-bounds}), we have that $w_{t}\in L^{\infty }\left(
\mathbb{R}_{\geq 0},H^{1/2}(\Gamma )\right) $. Thus, from the second
equation of (\ref{eq:pde-w}),
\begin{equation}
\Vert \partial _{\mathbf{n}}w\left( t\right) \Vert _{H^{1/2}(\Gamma )}\leq
\Vert w\left( t\right) \Vert _{H^{1/2}(\Gamma )}+\Vert w_{t}\left( t\right)
\Vert _{H^{1/2}(\Gamma )}\leq Q(R).  \label{eq:wt-bounds-3}
\end{equation}%
Combining (\ref{eq:wt-bounds-1}) and (\ref{eq:wt-bounds-3}) with (\ref%
{eq:H2-regularity}), and also applying (\ref{eq:wt-bounds}), proves that for
all $t\geq 0$,
\begin{equation*}
\Vert (w(t),w_{t}(t))\Vert _{\mathcal{D}_{\varepsilon }}\leq Q(R).
\end{equation*}

It follows that the operators $K_{\varepsilon }$ are uniformly compact (with
$t_{c}=0$).
\end{proof}

Next, we will discuss regularity properties of the weak solutions.

\begin{theorem}
\label{t:exponential-attraction-h} For each $\varepsilon \in (0,1]$, there
exists a closed and bounded subset $\mathcal{C}_{\varepsilon }\subset
\mathcal{D}_{\varepsilon } $, such that for every nonempty bounded subset $%
B\subset \mathcal{H}_{\varepsilon }$,%
\begin{equation}
\mathrm{dist}_{\mathcal{H}_{\varepsilon }}(S_{\varepsilon }(t)B,\mathcal{C}%
_{\varepsilon })\leq Q(\left\Vert B\right\Vert _{\mathcal{H}_{\varepsilon
}})e^{-\omega t},  \label{eq:transitivity}
\end{equation}%
where $Q$ and $\omega >0$ are independent of $\varepsilon .$
\end{theorem}

\begin{proof}
Let $\varepsilon \in (0,1]$. Define the subset $\mathcal{C}_{\varepsilon }$\
of $\mathcal{D}_{\varepsilon }$ by%
\begin{equation*}
\mathcal{C}_{\varepsilon }:=\left\{ \varphi \in \mathcal{D}_{\varepsilon
}:\Vert \varphi \Vert _{\mathcal{D}_{\varepsilon }}\leq Q(R)\right\} ,
\end{equation*}%
where $Q(R)>0$ is the function from Lemma \ref{t:uniform-compactness}, and $%
R>0$ is such that $\Vert \varphi _{0}\Vert _{\mathcal{H}_{\varepsilon }}\leq
R.$ Let now $\varphi _{0}=(u_{0},u_{1})\in \mathcal{B}_{\varepsilon }$
(endowed with the same topology of $\mathcal{H}_{\varepsilon }$). Then, for
all $t\geq 0$ and for all $\varphi _{0}\in \mathcal{B}_{\varepsilon }$, $%
S_{\varepsilon }(t)\varphi _{0}=Z_{\varepsilon }(t)\varphi
_{0}+K_{\varepsilon }(t)\varphi _{0}$, where $Z_{\varepsilon }(t)$ is
uniformly and exponentially decaying to zero by Lemma \ref{t:uniform-decay},
and, by Lemma \ref{t:uniform-compactness}, $K_{\varepsilon }(t)$ is
uniformly bounded in $\mathcal{D}_{\varepsilon }.$ In particular, there holds%
\begin{equation*}
\mathrm{dist}_{\mathcal{H}_{\varepsilon }}(S_{\varepsilon }(t)\mathcal{B}%
_{\varepsilon },\mathcal{C}_{\varepsilon })\leq Q(R)e^{-\omega t}.
\end{equation*}%
(Recall, $\omega >0$ is independent of $\varepsilon $ due to Lemma \ref%
{t:uniform-decay}).

Recall that, by Lemma \ref{t:bounded-absorbing-set-h}, we already know that
for each $\varepsilon \in (0,1]$ and for every nonempty bounded subset $B$
of $\mathcal{H}_{\varepsilon }$,
\begin{equation*}
\mathrm{dist}_{\mathcal{H}_{\varepsilon }}(S_{\varepsilon }(t)B,\mathcal{B}%
_{\varepsilon })\leq Q(R)e^{-\omega _{0}t},
\end{equation*}%
for all $t\geq 0$. In light of these estimates, (\ref{eq:transitivity}) can
now be accomplished by appealing to the transitivity property of the
exponential attraction (see, e.g., \cite[Theorem 5.1]{FGMZ04}). Note that (%
\ref{eq:transitivity}) entails that $\mathcal{C}_{\varepsilon }$ is a
compact attracting set in $\mathcal{H}_\varepsilon$ for $S_{\varepsilon }(t)$%
. The proof is finished.
\end{proof}

By standard arguments of the theory of attractors (see, e.g., \cite{Hale88,
Temam88}), the existence of a compact global attractor $\mathcal{A}%
_{\varepsilon }\subset \mathcal{C}_{\varepsilon }$ for the semigroup $%
S_{\varepsilon }(t)$ follows.

\begin{theorem}
\label{t:global-attractors-h} For each $\varepsilon \in (0,1]$, the semiflow
$S_{\varepsilon }$ generated by the solutions of the hyperbolic relaxation
problem (\ref{eq:pde-h-1})-(\ref{eq:pde-h-3}) admits a unique global
attractor%
\begin{equation*}
\mathcal{A}_{\varepsilon }=\omega (\mathcal{B}_{\varepsilon
}):=\bigcap_{s\geq t}{\overline{\bigcup_{t\geq 0}S_{\varepsilon }(t)\mathcal{%
B}_{\varepsilon }}}^{\mathcal{H}_{\varepsilon }}
\end{equation*}%
in $\mathcal{H}_{\varepsilon }$. Moreover, the following hold:

(i) For each $t\geq 0$, $S_{\varepsilon }(t)\mathcal{A}_{\varepsilon }=%
\mathcal{A}_{\varepsilon }$, and

(ii) For every nonempty bounded subset $B$ of $\mathcal{H}_{\varepsilon }$,%
\begin{equation}
\lim_{t\rightarrow \infty }\mathrm{dist}_{\mathcal{H}_{\varepsilon
}}(S_{\varepsilon }(t)B,\mathcal{A}_{\varepsilon })=0.  \label{gl_attraction}
\end{equation}

(iii) The global attractor $\mathcal{A}_{\varepsilon }$ is bounded in $%
\mathcal{D}_{\varepsilon }$ and trajectories on $\mathcal{A}_{\varepsilon }$
are strong solutions.
\end{theorem}

\begin{remark}
\label{nonlinear_bc}We can extend all the results in Sections \ref%
{well-posed}-\ref{global-hyper} (with the appropriate modifications, see
\cite{CGG11, GM12}) to the case when the linear boundary condition (\ref%
{eq:pde-h-2}) is replaced by%
\begin{equation*}
\partial _{\mathbf{n}}u+g\left( u\right) +u_{t}=0,\text{ on }\left( 0,\infty
\right) \times \Gamma ,
\end{equation*}%
such that $g\in C^{2}\left( \mathbb{R}\right) $ satisfies%
\begin{equation*}
|g^{^{\prime \prime }}\left( s\right) |\leq C_{g}\left( 1+\left\vert
s\right\vert ^{2}\right) ,\text{ }g^{^{\prime }}\left( s\right) \geq -\theta
_{g},\text{ }g\left( s\right) s\geq s^{2}-C_{g}^{^{\prime }},
\end{equation*}%
for all $s\in \mathbb{R}$, and some constants $C_{g}>0,$ $C_{g}^{^{\prime
}}\geq 0.$
\end{remark}

\subsection{The upper-semicontinuity of $\mathcal{A}_{\protect\varepsilon }$
for the singularly perturbed problem}

\label{s:continuity}

This section contains one of the main results of the paper, the proof of the
upper-semicontinuity of the family of global attractors given by the model
problems for $\varepsilon \in \lbrack 0,1]$. Recall that in the case $%
\varepsilon =0$, the limit parabolic problem, admits a global attractor $%
\mathcal{A}_{0}$ that is bounded in $\mathcal{V}^{2}$. Naturally, we will
study the continuity at $\varepsilon =0$. For $\varepsilon \in (0,1]$, we
know that the hyperbolic relaxation problem admits a global attractor $%
\mathcal{A}_{\varepsilon }$ in $\mathcal{D}_{\varepsilon }$. However, the
spaces involved with the parabolic problem invoke the trace of the solution
on the boundary $\Gamma $, whereas the spaces involved with the hyperbolic
relaxation problem do \emph{not} contain prescribed traces. Before we lift
the global attractor $\mathcal{A}_{0}$ for the parabolic problem into the
finite energy phase space for the hyperbolic relaxation problem, we need to
make an extension of $\mathcal{H}_{\varepsilon }$ so that it also includes
the information of the traces of $u$ and $u_{t}$.

To begin, we recall that the natural phase space for the parabolic problem (%
\ref{eq:pde-p-1})-(\ref{eq:pde-p-3}) is $Y=L^{2}(\Omega )\times L^{2}(\Gamma
),$ while the finite energy phase space for the hyperbolic relaxation
problem (\ref{eq:pde-h-1})-(\ref{eq:pde-h-3}) is $\mathcal{H}_{\varepsilon
}=H^{1}(\Omega )\times L^{2}(\Omega )$. Thus, we need to find a suitable
extension of the phase space for the hyperbolic relaxation problem so that,
when we lift the parabolic problem, both problems will be situated in the
\emph{same} framework. A natural way to make this extension is to introduce
the space
\begin{equation*}
\mathcal{X}_{0}=H^{1}(\Omega )\times L^{2}(\Gamma ),
\end{equation*}%
then the \emph{extended} phase space for the hyperbolic relaxation problem%
\begin{equation*}
\mathcal{X}_{\varepsilon }=\mathcal{X}_{0}\times Y=H^{1}(\Omega )\times
L^{2}(\Gamma )\times L^{2}(\Omega )\times L^{2}(\Gamma ).
\end{equation*}%
The space $\mathcal{X}_{\varepsilon }$ is Hilbert when endowed with the $%
\varepsilon $-weighted norm whose square is given by, for all $\zeta
=(u,\gamma ,v,\delta )\in \mathcal{X}_{\varepsilon }$,
\begin{equation*}
\Vert \zeta \Vert _{\mathcal{X}_{\varepsilon }}^{2}:=\Vert u\Vert
_{1}^{2}+\Vert \gamma \Vert _{L^{2}(\Gamma )}^{2}+\varepsilon \Vert v\Vert
^{2}+\varepsilon \Vert \delta \Vert _{L^{2}(\Gamma )}^{2}.
\end{equation*}%
It is then in the space $\mathcal{X}_{\varepsilon }$ that we can lift $%
\mathcal{A}_{0}$ and estimate the Hausdorff semidistance between (an
extension of) $\mathcal{A}_{\varepsilon }$ and $\mathcal{L}\mathcal{A}_{0}$
(for a proper lifting map $\mathcal{L}$) with the new extended topology.
However, it must be noted that the lifted attractor $\mathcal{L}\mathcal{A}%
_{0}$ is \emph{not} necessarily a global attractor when set in the extended
phase space. Finally, the topology that we will use to show the convergence
of the attractors at $\varepsilon =0$ will be defined with the
four-component norm of $\mathcal{X}_{\varepsilon }$.

For both problems, we also recall that trajectories on the attractor are
strong solutions due to the regularity results obtained in Sections \ref%
{s:parabolic} and \ref{global-hyper} (see Theorems \ref%
{t:global-attractor-and-regularity-p} and \ref{t:exponential-attraction-h}).
The regularized phase space $\mathcal{D}_{\varepsilon }$ for the hyperbolic
relaxation problem is isomorphically extended to%
\begin{align}
\widetilde{\mathcal{D}}_{\varepsilon }& :=\left\{ \left( u,\gamma ,v,\delta
\right) \in H^{2}\left( \Omega \right) \times H^{3/2}\left( \Gamma \right)
\times H^{1}\left( \Omega \right) \times H^{1/2}\left( \Gamma \right)
:\gamma =\mathrm{tr_{D}}\left( u\right) ,\right.   \label{extended_ph} \\
& \left. \delta =\mathrm{tr_{D}}\left( v\right) ,\text{ }\partial _{\mathbf{n%
}}u+\gamma =-\delta \text{ on }\Gamma \right\} .  \notag
\end{align}%
Of course, $\widetilde{\mathcal{D}}_{\varepsilon }\subset \mathcal{V}%
^{2}\times \mathcal{V}^{1}$ and the injection $\widetilde{\mathcal{D}}%
_{\varepsilon }\hookrightarrow \mathcal{X}_{\varepsilon }$ is compact.
Recall that, for each $\left( u_{0},u_{1}\right) \in \mathcal{D}%
_{\varepsilon }$, problem (\ref{eq:pde-h-1})-(\ref{eq:pde-h-3}) generates a
dynamical system $\left( S_{\varepsilon }\left( t\right) ,\mathcal{D}%
_{\varepsilon }\right) $ of strong solutions (cf., Theorem \ref%
{strong-sol-hyper}; see also \cite{Wu&Zheng06}). By appealing once more to
the continuity of the trace map $\mathrm{tr_{D}}:H^{s}\left( \Omega \right)
\rightarrow H^{s-1/2}\left( \Gamma \right) $, $s>1/2$, and exploiting the
results in Section \ref{well-posed}, it is not difficult to realize that we
can extend the semiflow $S_{\varepsilon }\left( t\right) $\ to a strongly
continuous semigroup%
\begin{equation}
\widetilde{S}_{\varepsilon }\left( t\right) :\widetilde{\mathcal{D}}%
_{\varepsilon }\rightarrow \widetilde{\mathcal{D}}_{\varepsilon },
\label{lifted-semigroup}
\end{equation}%
such that $\widetilde{S}_{\varepsilon }\left( t\right) $ is also Lipschitz
continuous in $\widetilde{\mathcal{D}}_{\varepsilon },$ endowed with the
metric topology of \ $\mathcal{V}^{2}\times \mathcal{V}^{1}$ (see Lemma \ref%
{Lipsch_higher} below). Recall that, by definition for $p,q\geq 1$,%
\begin{equation*}
\mathcal{V}^{p}\times \mathcal{V}^{q}=\left\{ \left( u,\gamma ,v,\delta
\right) \in H^{p}\left( \Omega \right) \times H^{p-1/2}\left( \Gamma \right)
\times H^{q}\left( \Omega \right) \times H^{q-1/2}\left( \Gamma \right)
:\gamma =\mathrm{tr_{D}}\left( u\right) ,\delta =\mathrm{tr_{D}}\left(
v\right) \right\} ,
\end{equation*}%
see Section \ref{s:parabolic} (as before, $\mathcal{V}^{p}\times \mathcal{V}%
^{q}$ is topologically isomorphic to $H^{p}\left( \Omega \right) \times
H^{q}\left( \Omega \right) $).

\begin{lemma}
\label{Lipsch_higher}Let $\varphi _{0},\theta _{0}\in \widetilde{\mathcal{D}}%
_{\varepsilon }$ such that $\left\Vert \varphi _{0}\right\Vert _{\widetilde{%
\mathcal{D}}_{\varepsilon }}\leq R,$ $\left\Vert \theta _{0}\right\Vert _{%
\widetilde{\mathcal{D}}_{\varepsilon }}\leq R,$ for every $\varepsilon \in
(0,1]$. Then the following estimate holds:%
\begin{equation}
\left\Vert \widetilde{S}_{\varepsilon }\left( t\right) \varphi _{0}-%
\widetilde{S}_{\varepsilon }\left( t\right) \theta _{0}\right\Vert _{%
\widetilde{\mathcal{D}}_{\varepsilon }}\leq \frac{Q\left( R\right) }{\sqrt{%
\varepsilon }}e^{\nu _{1}t}\left\Vert \varphi _{0}-\theta _{0}\right\Vert _{%
\widetilde{\mathcal{D}}_{\varepsilon }},  \label{lipschitzest2}
\end{equation}%
where $Q\left( R\right) >0$ and $\nu _{1}>0$ are independent of $\varepsilon
>0.$
\end{lemma}

\begin{proof}
Let $\varphi (t)=(u_{1}(t),u_{1\mid \Gamma }\left( t\right) ,\partial
_{t}u_{1}(t),\partial _{t}u_{1\mid \Gamma }(t))$ and $\theta
(t)=(u_{2}(t),u_{2\mid \Gamma }\left( t\right) ,\partial
_{t}u_{2}(t),\partial _{t}u_{2\mid \Gamma }(t))$ denote the corresponding
strong solutions with initial data $\varphi _{0}$ and $\theta _{0}$,
respectively. Then the difference $u\left( t\right) :=u_{1}\left( t\right)
-u_{2}\left( t\right) $ satisfies%
\begin{equation}
\left\{
\begin{array}{ll}
-\Delta u\left( t\right) =f^{^{\prime }}\left( u_{2}\left( t\right) \right)
-f^{^{\prime }}\left( u_{1}\left( t\right) \right) -u_{t}\left( t\right)
-\varepsilon u_{tt}\left( t\right) , & \text{a.e. in }\mathbb{R}_{+}\times
\Omega , \\
\partial _{\mathbf{n}}u\left( t\right) +u\left( t\right) =-u_{t}\left(
t\right) , & \text{a.e. in }\mathbb{R}_{+}\times \Gamma ,%
\end{array}%
\right.   \label{sys1}
\end{equation}%
subject to the initial condition%
\begin{equation*}
u\left( 0\right) =u_{1}\left( 0\right) -u_{2}\left( 0\right) .
\end{equation*}%
Setting $v:=\partial _{t}u_{1}-\partial _{t}u_{2}$, we have $\left(
v_{t},\psi \right) \in C^{1}\left( \left[ 0,T\right] \right) $ for every $%
\psi \in H^{1}\left( \Omega \right) $ (see the definition of strong
solution). Then $v$ solves the following identity%
\begin{align*}
& \frac{\mathrm{d}}{\mathrm{d}t}\left( \varepsilon v_{t}\left( t\right)
,\psi \right) +\left\langle \nabla v\left( t\right) ,\nabla \psi
\right\rangle +\left\langle v_{t}\left( t\right) ,\psi \right\rangle
+\left\langle v_{t}\left( t\right) +v\left( t\right) ,\psi \right\rangle
_{L^{2}\left( \Gamma \right) } \\
& =-\left\langle f^{^{\prime }}\left( u_{1}\left( t\right) \right)
-f^{^{\prime }}\left( u_{2}\left( t\right) \right) u_{1}\left( t\right)
,\psi \right\rangle -\left\langle f^{^{\prime }}\left( u_{2}\left( t\right)
\right) u\left( t\right) ,\psi \right\rangle ,
\end{align*}%
for almost all $t\in \left[ 0,T\right] .$ Testing with $\psi =v_{t}$, we
obtain%
\begin{align}
& \frac{1}{2}\frac{\mathrm{d}}{\mathrm{d}t}\left\{ \varepsilon \left\Vert
v_{t}\right\Vert ^{2}+\left\Vert \nabla v\right\Vert ^{2}+\left\Vert
v\right\Vert _{L^{2}\left( \Gamma \right) }^{2}\right\} +\left\Vert
v_{t}\right\Vert ^{2}+\left\Vert v_{t}\right\Vert _{L^{2}\left( \Gamma
\right) }^{2}  \label{lip2} \\
& =-\left\langle f^{^{\prime }}\left( u_{1}\right) -f^{^{\prime }}\left(
u_{2}\right) u_{1},v_{t}\right\rangle -\left\langle f^{^{\prime }}\left(
u_{2}\right) u,v_{t}\right\rangle .  \notag
\end{align}%
We can bound the terms on the right-hand side in a standard way,%
\begin{align*}
& \left\langle f^{^{\prime }}\left( u_{1}\right) -f^{^{\prime }}\left(
u_{2}\right) u_{1},v_{t}\right\rangle +\left\langle f^{^{\prime }}\left(
u_{2}\right) u,v_{t}\right\rangle  \\
& \leq Q\left( |u_{i}|_{\infty }\right) \left\Vert u\right\Vert ^{2}+\frac{1%
}{2}\left\Vert v_{t}\right\Vert ^{2}
\end{align*}%
(which follows easily on the account of the fact that $\left\Vert
(u_{i}\left( t\right) ,\partial _{t}u_{i}\left( t\right) )\right\Vert _{%
\mathcal{D}_{\varepsilon }}\leq R,$ $i=1,2,$ and the embedding $H^{2}\left(
\Omega \right) \hookrightarrow C^{0}\left( \overline{\Omega }\right) $),
then insert them into (\ref{lip2}). By virtue of (\ref%
{eq:continuous-dependence}) we get%
\begin{align}
& \varepsilon \left\Vert v_{t}\left( t\right) \right\Vert ^{2}+\left(
\left\Vert \nabla v\left( t\right) \right\Vert ^{2}+\left\Vert v\left(
t\right) \right\Vert _{L^{2}\left( \Gamma \right) }^{2}\right)   \label{lip3}
\\
& \leq Q\left( R\right) e^{\nu _{1}t}\left\Vert \varphi _{0}-\theta
_{0}\right\Vert _{\mathcal{X}_{\varepsilon }}^{2}+\varepsilon \left\Vert
v_{t}\left( 0\right) \right\Vert ^{2}+\left( \left\Vert \varphi _{0}-\theta
_{0}\right\Vert _{\widetilde{\mathcal{D}}_{\varepsilon }}^{2}\right) ,
\notag
\end{align}%
for almost all $t\in \left[ 0,T\right] .$ It remains to notice that, from (%
\ref{sys1}), there holds for every $\varepsilon \in (0,1],$%
\begin{align}
\varepsilon \left\Vert v_{t}\left( 0\right) \right\Vert ^{2}& \leq \frac{1}{%
\varepsilon }\left( \left\Vert \Delta u\left( 0\right) \right\Vert
^{2}+\left\Vert f^{^{\prime }}\left( u_{2}\left( 0\right) \right)
-f^{^{\prime }}\left( u_{1}\left( 0\right) \right) \right\Vert
^{2}+\left\Vert u_{t}\left( 0\right) \right\Vert ^{2}\right)   \label{lip4}
\\
& \leq \frac{Q\left( R\right) }{\varepsilon }\left\Vert \varphi _{0}-\theta
_{0}\right\Vert _{\widetilde{\mathcal{D}}_{\varepsilon }}^{2}  \notag
\end{align}%
Summing up, we obtain from (\ref{lip3})-(\ref{lip4}), that%
\begin{equation}
\varepsilon \left\Vert v_{t}\left( t\right) \right\Vert ^{2}+\left\Vert
v\left( t\right) \right\Vert _{1}^{2}\leq \frac{Q\left( R\right) }{%
\varepsilon }e^{\nu _{1}t}\left\Vert \varphi _{0}-\theta _{0}\right\Vert _{%
\widetilde{\mathcal{D}}_{\varepsilon }}^{2}.  \label{lip5}
\end{equation}%
We can now bound the term $\Vert u_{1}\left( t\right) -u_{2}\left( t\right)
\Vert _{2}$. As before, owing to standard elliptic regularity theory, we
have in (\ref{sys1}), using (\ref{lip5}), that%
\begin{eqnarray}
\Vert u\left( t\right) \Vert _{2}^{2} &\leq &C\left( \varepsilon
^{2}\left\Vert v_{t}\left( t\right) \right\Vert ^{2}+\left\Vert f^{^{\prime
}}\left( u_{2}\left( t\right) \right) -f^{^{\prime }}\left( u_{1}\left(
t\right) \right) \right\Vert ^{2}+\left\Vert v\left( t\right) \right\Vert
_{1}^{2}\right)   \label{lip6} \\
&\leq &Q\left( R\right) \left\Vert \varphi _{0}-\theta _{0}\right\Vert _{%
\widetilde{\mathcal{D}}_{\varepsilon }}^{2}.  \notag
\end{eqnarray}%
\newline
Finally, (\ref{lip5})-(\ref{lip6}) together with the fact that the trace map
\thinspace $H^{s}\left( \Omega \right) \rightarrow H^{s-1/2}\left( \Gamma
\right) $, $s>1/2$, is bounded yields the desired inequality (\ref%
{lipschitzest2}).
\end{proof}

By Lemma \ref{Lipsch_higher}, the family of global attractors $\left\{
\mathcal{A}_{\varepsilon }\right\} _{\varepsilon \in (0,1]}\subset \mathcal{D%
}_{\varepsilon }$ can be naturally extended to the family of compact sets $%
\left\{ \widetilde{\mathcal{A}}_{\varepsilon }\right\} _{\varepsilon \in
(0,1]},$%
\begin{equation}
\widetilde{\mathcal{A}}_{\varepsilon }=\left\{ \left( u,\gamma ,v,\delta
\right) \in \widetilde{\mathcal{D}}_{\varepsilon }:\left( u,v\right) \in
\mathcal{A}_{\varepsilon }\right\}   \label{conj_gl}
\end{equation}%
which are bounded in $\widetilde{\mathcal{D}}_{\varepsilon }$ and compact in
$\mathcal{X}_{\varepsilon }$. Note that we do not claim that $\widetilde{%
\mathcal{A}}_{\varepsilon }$ is a global attractor for $(\widetilde{S}%
_{\varepsilon }\left( t\right) ,\mathcal{X}_{\varepsilon })$, see Remark \ref%
{rem_notglob} below. Also, it is in the space $\mathcal{V}^{2}\times
Y\subset \mathcal{X}_{\varepsilon }$ where we lift the parabolic problem.
Since the global attractor $\mathcal{A}_{0}$ for (\ref{eq:pde-p-1})-(\ref%
{eq:pde-p-3}) is a bounded subset of the space $\mathcal{V}^{2}\subset
C\left( \overline{\Omega }\right) \times C\left( \Gamma \right) $ (since $%
\Omega \subset \mathbb{R}^{3}$), the canonical extension map%
\begin{equation}
\mathcal{E}:\mathcal{V}^{2}\rightarrow Y
\label{eq:canonical-extension-map-1}
\end{equation}%
is well-defined with
\begin{equation}
(u,u_{\mid \Gamma })\mapsto (\Delta u-f(u),-\partial _{\mathbf{n}}u-u_{\mid
\Gamma }),  \label{eq:canonical-extension-map-2}
\end{equation}%
and so the corresponding lift map%
\begin{equation}
\mathcal{L}:\mathcal{V}^{2}\rightarrow \mathcal{V}^{2}\times Y
\label{eq:lift-map-1}
\end{equation}%
is defined by%
\begin{equation}
(u,u_{\mid \Gamma })\mapsto (u,u_{\mid \Gamma },\Delta u-f(u),-\partial _{%
\mathbf{n}}u-u_{\mid \Gamma }).  \label{eq:lift-map-2}
\end{equation}

Let $\mathcal{A}_{0}$ denote the global attractor of the limit parabolic
problem (see Theorem \ref{t:global-attractor-and-regularity-p}) and let $%
\widetilde{\mathcal{A}}_\varepsilon$, $\varepsilon\in(0,1]$, be the sets
defined in (\ref{conj_gl}). Define the family of compact sets in $\mathcal{X}%
_{\varepsilon }$ by%
\begin{equation}
\mathbb{A}_{\varepsilon }:=\left\{
\begin{array}{ll}
\widetilde{\mathcal{A}}_{0}:=\mathcal{L}\mathcal{A}_{0} & \text{for}%
~\varepsilon =0 \\
\widetilde{\mathcal{A}}_{\varepsilon } & \text{for}~\varepsilon \in (0,1].%
\end{array}%
\right.  \label{eq:family}
\end{equation}

\begin{remark}
\label{rem_notglob}The compact set $\widetilde{\mathcal{A}}_{\varepsilon }$\
is \emph{not} a global attractor for $\widetilde{S}_{\varepsilon }\left(
t\right) $ acting on the phase-space $\mathcal{X}_{\varepsilon }$ since
traces of functions in $L^{2}\left( \Omega \right) $ are not well-defined in
$L^{2}\left( \Gamma \right) $. By construction (\ref{conj_gl}), $\widetilde{%
\mathcal{A}}_{\varepsilon }$ is only topologically conjugated to the global
attractor $\mathcal{A}_{\varepsilon }$ associated with the dynamical system $%
\left( S_{\varepsilon },\mathcal{H}_{\varepsilon }\right) .$
\end{remark}

The main result of this section can be now stated as follows.

\begin{theorem}
\label{t:upper-continuity} The family $\left\{\mathbb{A}_{\varepsilon
}\right\}_{\varepsilon \in \lbrack 0,1]},$ defined by (\ref{eq:family}), is
upper-semicontinuous at $\varepsilon =0$ in the topology of $\mathcal{X}_{1}$%
. More precisely, there holds
\begin{equation}
\lim_{\varepsilon \rightarrow 0}\mathrm{dist}_{\mathcal{X}_{1}}(\mathbb{A}%
_{\varepsilon },\mathbb{A}_{0}):=\lim_{\varepsilon \rightarrow 0}\sup_{a\in
\widetilde{\mathcal{A}}_{\varepsilon }}\inf_{b\in \widetilde{\mathcal{A}}%
_{0}}\Vert a-b\Vert _{\mathcal{X}_{1}}=0.  \label{eq:upper-semicontinuity}
\end{equation}
\end{theorem}

\begin{proof}
Our proof essentially follows the classical argument in \cite%
{Hale&Raugel88,Hale88} and also \cite[Theorem 3.31]{Milani&Koksch05}. Of
course, modifications are required to account for the terms on the boundary.
Let $\zeta =(u,\gamma ,v,\delta )\in \widetilde{\mathcal{A}}_{\varepsilon }$
and $\bar{\zeta}=(\bar{u},\bar{\gamma},\bar{v},\bar{\delta})\in \widetilde{%
\mathcal{A}}_{0}$. We need to show that
\begin{equation}
\begin{aligned} \label{eq:upper-1}
\sup_{(u,\gamma,v,\delta)\in\widetilde{\mathcal{A}}_\varepsilon}\inf_{(\bar
u,\bar \gamma,\bar v,\bar \delta)\in\widetilde{\mathcal{A}}_0} & \left(
\|u-\bar u\|^2_1 + \|\gamma-\bar \gamma\|^2_{L^2(\Gamma)} + \right. \\ &
\left. + \|v-\bar v\|^2 + \|\delta-\bar \delta\|^2_{L^2(\Gamma)}
\right)^{1/2} \rightarrow 0 ~\text{as}~ \varepsilon\rightarrow 0.
\end{aligned}
\end{equation}%
Assuming to the contrary that (\ref{eq:upper-1}) did not hold, then there
exist $\eta _{0}>0$ and sequences $(\varepsilon _{n})_{n\in \mathbb{N}%
}\subset (0,1]$, $(\zeta _{n})_{n\in \mathbb{N}}=((u_{n},\gamma
_{n},v_{n},\delta _{n}))_{n\in \mathbb{N}}\subset \widetilde{\mathcal{A}}%
_{\varepsilon _{n}}$, such that $\varepsilon _{n}\rightarrow 0$ and for all $%
n\in \mathbb{N}$,
\begin{equation}
\inf_{(\bar{u},\bar{\gamma},\bar{v},\bar{\delta})\in \widetilde{\mathcal{A}}%
_{0}}\left( \Vert u_{n}-\bar{u}\Vert _{1}^{2}+\Vert \gamma _{n}-\bar{\gamma}%
\Vert _{L^{2}(\Gamma )}^{2}+\Vert v_{n}-\bar{v}\Vert ^{2}+\Vert \delta _{n}-%
\bar{\delta}\Vert _{L^{2}(\Gamma )}^{2}\right) \geq \eta _{0}^{2}.
\label{eq:upper-2}
\end{equation}%
By Theorem \ref{t:exponential-attraction-h}, the compact sets $\widetilde{%
\mathcal{A}}_{\varepsilon _{n}}$ are bounded in the space $\widetilde{%
\mathcal{D}}_{1}$ (see (\ref{extended_ph}) with $\varepsilon=1$) and we have
the following uniform bound, for some positive constant $C>0$ independent of
$n$,
\begin{equation*}
\Vert u_{n}\Vert _{2}^{2}+\Vert \gamma _{n}\Vert _{H^{3/2}(\Gamma
)}^{2}+\Vert v_{n}\Vert _{1}^{2}+\Vert \delta _{n}\Vert _{H^{1/2}(\Gamma
)}^{2}\leq C.
\end{equation*}%
This means that there is a weakly converging subsequence of $(\zeta
_{n})_{n\in \mathbb{N}}$ (not relabelled) that converges to some $(u^{\ast
},\gamma ^{\ast },v^{\ast },\delta ^{\ast })$ weakly in $\widetilde{\mathcal{%
D}}_{1}$. By the compactness of the embedding $\widetilde{\mathcal{D}}%
_{1}\hookrightarrow \mathcal{X}_{1}$, the subsequence converges strongly in $%
\mathcal{X}_{1}$. It now suffices to show that $(u^{\ast },\gamma ^{\ast
},v^{\ast },\delta ^{\ast })\in \widetilde{\mathcal{A}}_{0},$ since this is
a contradiction to (\ref{eq:upper-2}).

With each $\zeta _{n}=(u_{n},\gamma _{n},v_{n},\delta _{n})\in \widetilde{%
\mathcal{A}}_{\varepsilon _{n}}$, then, for each $n\in \mathbb{N}$, there is
a complete orbit
\begin{equation*}
(u^{n}(t),u^{n}_{\mid \Gamma }(t),u_{t}^{n}(t),u^{n}_{t\mid \Gamma
}(t))_{t\in \mathbb{R}}=(\widetilde{S}_{\varepsilon _{n}}(t)(u_{n},\gamma
_{n},v_{n},\delta _{n}))_{t\in \mathbb{R}}
\end{equation*}%
contained in $\widetilde{\mathcal{A}}_{\varepsilon _{n}}$ and passing
through $(u_{n},\gamma _{n},v_{n},\delta _{n})$ where
\begin{equation*}
(u^{n}(0),u^{n}_{\mid \Gamma }(0),u_{t}^{n}(0),u^{n}_{t\mid \Gamma
}(0))=(u_{n},\gamma _{n},v_{n},\delta _{n})
\end{equation*}%
(cf., e.g., \cite[Proposition 2.39]{Milani&Koksch05}).

In view of the regularity $\widetilde{\mathcal{A}}_{\varepsilon _{n}}\subset
\widetilde{\mathcal{D}}_{1}$ (see (\ref{eq:regularity-property})), we obtain
the uniform bounds:%
\begin{equation}
\varepsilon _{n}\Vert u_{tt}^{n}(t)\Vert ^{2}+\Vert u_{t}^{n}(t)\Vert
_{1}^{2}+\Vert u_{t}^{n}(t)\Vert _{H^{1/2}(\Gamma )}^{2}+\Vert u^{n}(t)\Vert
_{2}^{2}+\Vert u^{n}(t)\Vert _{H^{3/2}(\Gamma )}^{2}\leq C,
\label{eq:upper-3}
\end{equation}%
where the constant $C>0$ is independent of $t$ and $\varepsilon _{n}$. Now,
for all $T>0$, the functions $u^{\varepsilon _{n}}$, $u_{\mid \Gamma
}^{\varepsilon _{n}}$, $u_{t}^{\varepsilon _{n}}$, $u_{t\mid \Gamma
}^{\varepsilon _{n}}$ and $\sqrt{\varepsilon _{n}}u_{tt}^{\varepsilon _{n}}$
are, respectively, bounded in $L^{\infty }(-T,T;H^{2}(\Omega ))$, $L^{\infty
}(-T,T;H^{3/2}(\Gamma ))$, $L^{\infty }(-T,T;H^{1}(\Omega ))$, $L^{\infty
}(-T,T;H^{1/2}(\Gamma ))$ and $L^{\infty }(-T,T;L^{2}(\Omega ))$. Thus,
there is a function $u$ and a subsequence (not relabelled), in which,
\begin{equation}
u^{\varepsilon _{n}}\rightharpoonup u~\text{in}~L^{\infty
}(-T,T;H^{2}(\Omega ))~(\text{weakly*}),  \label{eq:conv-1}
\end{equation}%
\begin{equation}
u_{\mid \Gamma }^{\varepsilon _{n}}\rightharpoonup u_{\mid \Gamma }~\text{in}%
~L^{\infty }(-T,T;H^{3/2}(\Gamma ))~(\text{weakly*}),  \label{eq:conv-2}
\end{equation}%
\begin{equation}
u_{t}^{\varepsilon _{n}}\rightharpoonup u_{t}~\text{in}~L^{\infty
}(-T,T;H^{1}(\Omega ))~(\text{weakly*}),  \label{eq:conv-3}
\end{equation}%
\begin{equation}
u_{t\mid \Gamma }^{\varepsilon _{n}}\rightharpoonup u_{t\mid \Gamma }~\text{%
in}~L^{\infty }(-T,T;H^{1/2}(\Gamma ))~(\text{weakly*}),  \label{eq:conv-4}
\end{equation}%
\begin{equation}
\varepsilon _{n}u_{tt}^{\varepsilon _{n}}\rightarrow 0~\text{in}~L^{\infty
}(-T,T;L^{2}(\Omega ))~(\text{strongly}).  \label{eq:conv-5}
\end{equation}

The above convergence properties yield%
\begin{equation}
u^{\varepsilon _{n}}\rightarrow u~\text{in}~C(-T,T;H^{1}(\Omega ))~(\text{%
strongly})  \label{eq:conv-6}
\end{equation}%
owing to the following embedding%
\begin{equation}
\{u\in L^{\infty }(-T,T;H^{2}(\Omega )):u_{t}\in L^{\infty
}(-T,T;H^{1}(\Omega ))\}\hookrightarrow C(-T,T;H^{2-\eta }(\Omega )),
\label{eq:injection-1}
\end{equation}%
which is compact for every $\eta \in \left( 0,1\right) $ (see, e.g., \cite%
{Lions69}). The strong property (\ref{eq:conv-6}) allows us to identify the
correct limit in the nonlinear term when $\varepsilon _{n}\rightarrow 0.$
Moreover, from (\ref{eq:conv-1}) and the fact that $H^{2}\left( \Omega
\right) \hookrightarrow C^{0}\left( \overline{\Omega }\right) $, it follows
that%
\begin{align}
\sup_{t\in \left[ -T,T\right] }\left\Vert f\left( u^{\varepsilon
_{n}}\right) -f\left( u\right) \right\Vert ^{2}& \leq \sup_{t\in \left[ -T,T%
\right] }Q_{\ast }\left( \left\vert u^{\varepsilon _{n}}\left( t\right)
\right\vert _{\infty },\left\vert u\left( t\right) \right\vert _{\infty
}\right) \left\Vert u^{\varepsilon _{n}}\left( t\right) -u\left( t\right)
\right\Vert ^{2} \\
& \leq C\left( \Omega \right) \sup_{t\in \left[ -T,T\right] }\left\Vert
u^{\varepsilon _{n}}\left( t\right) -u\left( t\right) \right\Vert ^{2},
\notag
\end{align}%
for some positive (increasing) function $Q_{\ast }:\mathbb{R}_{+}\times
\mathbb{R}_{+}\rightarrow \mathbb{R}_{+}$, independent of $n$ and $%
\varepsilon _{n}.$ By virtue of (\ref{eq:conv-6}) it is then easy to see that%
\begin{equation*}
f(u^{\varepsilon _{n}})\rightarrow f(u)~\text{in}~C(-T,T;L^{2}\left( \Omega
\right) )~(\text{strongly}).
\end{equation*}

It follows that $u$ is a weak solution of the limit parabolic problem on $%
\mathbb{R}$. In particular, $(u_{n},\gamma _{n})=(u^{n}(0),u^{n}_{\mid
\Gamma }(0))\rightarrow (u(0),u_{\mid \Gamma }(0))$ in $\mathcal{V}^{1}$.
Hence, we have that $(u(0),u_{\mid \Gamma }(0))=(u^{\ast },\gamma ^{\ast })$%
, and therefore $(u(0),u_{\mid \Gamma }(0))\in \mathcal{V}^{2}$. As $%
(u,u_{\mid\Gamma})$ is a complete orbit through $(u^{\ast },\gamma ^{\ast })$%
, it follows that $(u^{\ast },\gamma ^{\ast })\in \mathcal{A}_{0}$. It
remains to show that $v^{\ast }=\Delta u^{\ast }-f(u^{\ast })$ and $\delta
^{\ast }=-\partial _{\mathbf{n}}\gamma ^{\ast }-\gamma ^{\ast }$, in which
case $(u^{\ast },\gamma ^{\ast },v^{\ast },\delta ^{\ast })\in \widetilde{%
\mathcal{A}}_{0}$.

Now by (\ref{eq:conv-5}) and (\ref{eq:upper-3}), it follows that
\begin{equation*}
\|\varepsilon_n u^n_{tt}(0)\|=\sqrt{\varepsilon_n}\|\sqrt{\varepsilon_{n}}%
u_{tt}^{n}(0)\| \leq \sqrt{\varepsilon _{n}}C,
\end{equation*}
and so $\varepsilon _{n}u_{tt}^{n}(0)\rightarrow 0$ in $L^{2}(\Omega )$ as $%
\varepsilon _{n}\rightarrow 0$. With this at hand,
\begin{equation*}
\begin{aligned} u^n_t(0) & = - \varepsilon_n u^n_{tt}(0) + \Delta u^n(0) -
f(u^n(0)) \\ & = - \varepsilon_n u^n_{tt}(0) + \Delta u^* - f(u^*),
\end{aligned}
\end{equation*}%
so that
\begin{equation}
u_{t}^{n}(0)\rightharpoonup \Delta u^{\ast }-f(u^{\ast })~\text{in}%
~L^{2}(\Omega )~(\text{weakly}).  \label{eq:upper-4}
\end{equation}%
Since $u_{t}^{n}(0)=v^{n}$, then with (\ref{eq:upper-4}) we have that
\begin{equation}
v^{\ast }=\Delta u^{\ast }-f(u^{\ast }).  \label{eq:upper-5}
\end{equation}%
Similarly, since
\begin{equation*}
u_{t\mid \Gamma }^{n}(0)=-\partial _{\mathbf{n}}u^{n}(0)-u_{\mid \Gamma
}^{n}(0)
\end{equation*}%
and since $u_{\mid \Gamma }^{n}(0)=\gamma ^{\ast }$ and $u_{t\mid \Gamma
}^{n}(0)=\delta ^{\ast }$, then
\begin{equation}
\delta ^{\ast }=-\partial _{\mathbf{n}}\gamma ^{\ast }-\gamma ^{\ast }.
\label{eq:upper-6}
\end{equation}%
We know $(u^{\ast },\gamma ^{\ast })\in \mathcal{A}_{0}$, so (\ref%
{eq:upper-5}) and (\ref{eq:upper-6}) imply that $(u^{\ast },\gamma ^{\ast
},v^{\ast },\delta ^{\ast })\in \widetilde{\mathcal{A}}_{0},$ in
contradiction to (\ref{eq:upper-2}). This proves the assertion and completes
the proof.
\end{proof}

\section{Exponential attractors}

\label{s:exponential-attractors}

Exponential attractors (sometimes called, inertial sets) are positively
invariant sets possessing finite fractal dimension that attract bounded
subsets of the phase space exponentially fast. It can readily be seen that
when both a global attractor $\mathcal{A} $ and an exponential attractor $%
\mathcal{M}$ exist, then $\mathcal{A}\subseteq \mathcal{M}$, and so the
global attractor is also finite dimensional. The existence of an exponential
attractor depends on certain properties of the semigroup; namely, the
smoothing property for the difference of any two trajectories and the
existence of a more regular bounded absorbing set in the phase space (see,
e.g., \cite{EFNT95}, \cite{EMZ00}).

The main result of this section is the following.

\begin{theorem}
\label{t:exponential-attractors-h} For each $\varepsilon \in (0,1]$, the
dynamical system $\left( S_{\varepsilon },\mathcal{H}_{\varepsilon }\right) $
associated with (\ref{eq:pde-h-1})-(\ref{eq:pde-h-3}) admits an exponential
attractor $\mathcal{M}_{\varepsilon }$ compact in $\mathcal{H}_{\varepsilon
},$ and bounded in $\mathcal{C}_{\varepsilon }.$ Moreover, there hold:

(i) For each $t\geq 0$, $S_{\varepsilon }(t)\mathcal{M}_{\varepsilon
}\subseteq \mathcal{M}_{\varepsilon }$.

(ii) The fractal dimension of $\mathcal{M}_{\varepsilon }$ with respect to
the metric $\mathcal{H}_{\varepsilon }$ is finite, namely,%
\begin{equation*}
\dim _{F}\left( \mathcal{M}_{\varepsilon },\mathcal{H}_{\varepsilon }\right)
\leq C_{\varepsilon }<\infty ,
\end{equation*}%
for some positive constant $C_{\varepsilon }$ which \emph{depends} on $%
\varepsilon .$

(iii) There exist $\varrho >0$ and a positive nondecreasing function $%
Q_{\varepsilon }$ such that, for all $t\geq 0$,
\begin{equation*}
\mathrm{dist}_{\mathcal{H}_{\varepsilon }}(S_{\varepsilon }(t)B,\mathcal{M}%
_{\varepsilon })\leq Q_{\varepsilon }(\Vert B\Vert _{\mathcal{H}%
_{\varepsilon }})e^{-\varrho t},
\end{equation*}%
for every nonempty bounded subset $B$ of $\mathcal{H}_{\varepsilon }.$
\end{theorem}

\begin{remark}
Above,%
\begin{equation*}
\dim _{\mathrm{F}}(\mathcal{M}_{\varepsilon },\mathcal{H}_{\varepsilon
}):=\limsup_{r\rightarrow 0}\frac{\ln \mu _{\mathcal{H}_{\varepsilon }}(%
\mathcal{M}_{\varepsilon },r)}{-\ln r}<\infty ,
\end{equation*}%
where, $\mu _{\mathcal{H}_{\varepsilon }}(\mathcal{Z},r)$ denotes the
minimum number of $r$-balls from $\mathcal{H}_{\varepsilon }$ required to
cover $\mathcal{Z}\subset \mathcal{H}_{\varepsilon }$.
\end{remark}

\begin{corollary}
There holds%
\begin{equation*}
\dim _{\mathrm{F}}(\mathcal{A}_{\varepsilon },\mathcal{H}_{\varepsilon
})\leq \dim _{\mathrm{F}}(\mathcal{M}_{\varepsilon },\mathcal{H}%
_{\varepsilon }).
\end{equation*}%
As a consequence, $\mathcal{A}_{\varepsilon }$ has finite fractal dimension
which depends on $\varepsilon >0$.
\end{corollary}

\begin{remark}
Unfortunately, we cannot show that the fractal dimension of $\mathcal{M}%
_{\varepsilon }$ is uniform with respect to $\varepsilon >0$ (see the
subsequent lemmas).
\end{remark}

The proof of Theorem \ref{t:exponential-attractors-h} follows from the
application of an abstract result tailored specifically to our needs (see,
e.g., \cite[Proposition 1]{EMZ00}, \cite{FGMZ04}, \cite{GGMP05}; cf. also
Remark \ref{rem_att}\ below).

\begin{proposition}
\label{abstract1}Let $\left( S_{\varepsilon },\mathcal{H}_{\varepsilon
}\right) $ be a dynamical system for each $\varepsilon >0$. Assume the
following hypotheses hold:

\begin{enumerate}
\item[(C1)] There exists a bounded absorbing set $\mathcal{B}_{\varepsilon
}^{1}\subset \mathcal{D}_{\varepsilon }$ which is positively invariant for $%
S_{\varepsilon }\left( t\right) .$ More precisely, there exists a time $%
t_{1}>0,$ which \emph{depends} on $\varepsilon >0$, such that%
\begin{equation*}
S_{\varepsilon }(t)\mathcal{B}_{\varepsilon }^{1}\subset \mathcal{B}%
_{\varepsilon }^{1}
\end{equation*}%
for all $t\geq t_{1}$ where $\mathcal{B}_{\varepsilon }^{1}$ is endowed with
the topology of $\mathcal{H}_{\varepsilon }.$

\item[(C2)] There is $t^{\ast }\geq t_{1}$ such that the map $S_{\varepsilon
}(t^{\ast })$ admits the decomposition, for each $\varepsilon \in (0,1]$ and
for all $\varphi _{0},\theta _{0}\in \mathcal{B}_{\varepsilon }^{1}$,
\begin{equation*}
S_{\varepsilon }(t^{\ast })\varphi _{0}-S_{\varepsilon }(t^{\ast })\theta
_{0}=L_{\varepsilon }(\varphi _{0},\theta _{0})+R_{\varepsilon }(\varphi
_{0},\theta _{0})
\end{equation*}%
where, for some constants $\alpha ^{\ast }\in (0,\frac{1}{2})$ and $\Lambda
^{\ast }=\Lambda ^{\ast }(\Omega ,t^{\ast })\geq 0$ with $\Lambda ^{\ast }$
depending on $\varepsilon >0$, the following hold:%
\begin{equation}
\Vert L_{\varepsilon }(\varphi _{0},\theta _{0})\Vert _{\mathcal{H}%
_{\varepsilon }}\leq \alpha ^{\ast }\Vert \varphi _{0}-\theta _{0}\Vert _{%
\mathcal{H}_{\varepsilon }}  \label{eq:difference-decomposition-L}
\end{equation}%
and%
\begin{equation}
\Vert R_{\varepsilon }(\varphi _{0},\theta _{0})\Vert _{\mathcal{D}%
_{\varepsilon }}\leq \Lambda ^{\ast }\Vert \varphi _{0}-\theta _{0}\Vert _{%
\mathcal{H}_{\varepsilon }}.  \label{eq:difference-decomposition-K}
\end{equation}

\item[(C3)] The map%
\begin{equation*}
(t,U)\mapsto S_{\varepsilon }(t)U:[t^{\ast },2t^{\ast }]\times \mathcal{B}%
_{\varepsilon }^{1}\rightarrow \mathcal{B}_{\varepsilon }^{1}
\end{equation*}%
is Lipschitz continuous on $\mathcal{B}_{\varepsilon }^{1}$ in the topology
of $\mathcal{H}_{\varepsilon }$.
\end{enumerate}

Then, $\left( S_{\varepsilon },\mathcal{H}_{\varepsilon }\right) $ possesses
an exponential attractor $\mathcal{M}_{\varepsilon }$ in $\mathcal{B}%
_{\varepsilon }^{1}.$
\end{proposition}

We now show that the assumptions (C1)-(C3) hold for $\left( S_{\varepsilon
}\left( t\right) ,\mathcal{H}_{\varepsilon }\right) $. We begin with a
higher-order dissipative estimate in the norm of $\mathcal{D}_{\varepsilon
}. $

\begin{lemma}
\label{t:to-H1} Condition (C1) holds for fixed $\varepsilon \in (0,1]$.
\end{lemma}

\begin{proof}
The proof is very similar to the proof of Lemma \ref{t:uniform-compactness}.
Indeed, let $\varepsilon \in (0,1]$, $\varphi _{0}=(u_{0},u_{1})\in \mathcal{%
D}_{\varepsilon }$ and $\varphi \left( t\right) =S_{\varepsilon }\left(
t\right) \varphi _{0}$. In this setting, we differentiate (\ref{eq:pde-h-1}%
)-(\ref{eq:pde-h-3}) with respect to $t$ and let $h=u_{t}$. We set $\beta $
in (\ref{eq:beta}) to be $\beta =\vartheta $ where we recall that $\vartheta
>0$ is due to assumption (\ref{eq:assumption-f-3}). Then we easily obtain
the analogue of the differential inequality (\ref{eq:for-C1}) except that
the size of the initial data now depends on the norm of $\mathcal{D}%
_{\varepsilon }$, i.e., $\varphi _{0}=(u_{0},u_{1})\in \mathcal{D}%
_{\varepsilon }$ (here the initial conditions are not necessarily equal to
zero). Thus, after applying (\ref{eq:diff-lower}) and (\ref{eq:diff-upper}),
there exist a positive and nondecreasing function $Q$ and a constant $C>0$
such that%
\begin{equation}
\left\Vert \left( h(t),h_{t}(t)\right) \right\Vert _{\mathcal{H}%
_{\varepsilon }}^{2}\leq Q(\Vert \left( h(0),h_{t}(0)\right) \Vert _{%
\mathcal{H}_{\varepsilon }})e^{-\omega _{6}t/2}+C\left( R\right)
\label{eq:uniform-reg-bound}
\end{equation}%
($Q$, $\omega _{6}$ and $C$ are independent of $\varepsilon $) with $R>0$
such that $\Vert \varphi _{0}\Vert _{\mathcal{H}_{\varepsilon }}\leq R$.
Arguing as in Theorem \ref{t:uniform-compactness} by exploiting $H^{2}$%
-elliptic regularity theory, we also deduce%
\begin{equation}
\left\Vert \varphi \left( t\right) \right\Vert _{\mathcal{D}_{\varepsilon
}}^{2}\leq Q_{\varepsilon }(\Vert \varphi _{0}\Vert _{\mathcal{D}%
_{\varepsilon }})e^{-\omega _{6}t/2}+C\left( R\right) ,  \label{unif_h2}
\end{equation}%
for some new function $Q_{\varepsilon }$ which depends on $\varepsilon >0$.
Indeed, using the equations (\ref{eq:pde-h-1})-(\ref{eq:pde-h-3}), it is not
difficult to show that there holds $\Vert \left( h(0),h_{t}(0)\right) \Vert
_{\mathcal{H}_{\varepsilon }}\leq \frac{C}{\sqrt{\varepsilon }}\Vert \varphi
_{0}\Vert _{\mathcal{D}_{\varepsilon }}$, whence (\ref{unif_h2}).
Consequently, there exists $R_{1}>0$ (independent of time, $\varepsilon >0$
and initial data) such that $S_{\varepsilon }\left( t\right) $ possesses an
absorbing ball $\mathcal{B}_{\varepsilon }^{1}=B_{\mathcal{D}_{\varepsilon
}}\left( R_{1}\right) $ of radius $R_{1}$ centered at $0$, which is bounded
in $\mathcal{D}_{\varepsilon }$. This establishes condition (C1).
\end{proof}

\begin{remark}
\label{r:reg-exp-attr} \label{expo_att_set}Unfortunately, the bound in the
space $\mathcal{D}_{\varepsilon }$ is not uniform as $\varepsilon
\rightarrow 0^{+}$. Indeed the function $Q_{\varepsilon }\left( \cdot
\right) $ in (\ref{unif_h2}) blows up as $\varepsilon \rightarrow 0^{+}$.
Finally, arguing in a standard way as in Theorem \ref%
{t:exponential-attraction-h}, $\mathcal{B}_{\varepsilon }^{1}$ is in fact
exponentially attracting in $\mathcal{H}_{\varepsilon }$.
\end{remark}

\begin{lemma}
\label{t:to-C2} Condition (C2) holds for each fixed $\varepsilon \in (0,1]$.
\end{lemma}

\begin{proof}
Let $\varepsilon \in (0,1]$. Let $\varphi _{0},\theta _{0}\in \mathcal{B}%
_{\varepsilon }^{1}$. Define the pair of trajectories, for $t\geq 0$, $%
\varphi (t)=S_{\varepsilon }(t)\varphi _{0}=(u(t),u_{t}(t))$ and $\theta
(t)=S_{\varepsilon }(t)\theta _{0}=(v(t),v_{t}(t))$. For each $t\geq 0$,
decompose the difference $\bar{\zeta}(t):=\varphi (t)-\theta (t)$ with $\bar{%
\zeta}_{0}:=\varphi _{0}-\theta _{0}$ as follows:%
\begin{equation*}
\bar{\zeta}(t)=\bar{\varphi}(t)+\bar{\theta}(t)
\end{equation*}%
where $\bar{\varphi}(t)=(\bar{u}(t),\bar{u}_{t}(t))$ and $\bar{\theta}(t)=(%
\bar{v}(t),\bar{v}_{t}(t))$ are solutions of the problems:%
\begin{equation}
\left\{
\begin{array}{ll}
\varepsilon \bar{u}_{tt}+\bar{u}_{t}-\Delta \bar{u}=0 & \text{in}~(0,\infty
)\times \Omega \\
\partial _{\mathbf{n}}\bar{u}+\bar{u}+\bar{u}_{t}=0 & \text{on}~(0,\infty
)\times \Gamma \\
\bar{\varphi}(0)=\varphi _{0}-\theta _{0} & \text{in}~\Omega%
\end{array}%
\right.  \label{eq:diff-decomp-u}
\end{equation}%
and
\begin{equation}
\left\{
\begin{array}{ll}
\varepsilon \bar{v}_{tt}+\bar{v}_{t}-\Delta \bar{v}=f(v)-f(u) & \text{in}%
~(0,\infty )\times \Omega \\
\partial _{\mathbf{n}}\bar{v}+\bar{v}+\bar{v}_{t}=0 & \text{on}~(0,\infty
)\times \Gamma \\
\bar{\theta}(0)=0 & \text{in}~\Omega .%
\end{array}%
\right.  \label{eq:diff-decomp-v}
\end{equation}

By estimating along the usual lines, multiplying (\ref{eq:diff-decomp-u})$_1$
by $2\bar{u}_{t}+\bar{u}$ in $L^{2}(\Omega )$, we easily obtain the
differential inequality, for almost all $t\geq 0$,
\begin{equation}
\frac{\mathrm{d}}{\mathrm{d}t}\mathcal{N}_{\varepsilon }+\omega_7\mathcal{N}%
_{\varepsilon }\leq 0,  \label{eq:diff-decomp-1}
\end{equation}%
for some positive constant $\omega_7$ sufficiently small and independent of $%
\varepsilon$, and for
\begin{equation}  \label{eq:H-equivalent-norm}
\mathcal{N}_{\varepsilon }=\mathcal{N}_{\varepsilon }(\bar{\varphi}%
(t)):=\varepsilon \Vert \bar{u}_{t}(t)\Vert ^{2}+\varepsilon\langle \bar{u}%
_{t}(t),\bar{u}(t)\rangle +\Vert \nabla \bar{u}(t)\Vert ^{2}+\Vert \bar{u}%
(t)\Vert _{L^{2}(\Gamma )}^{2}.
\end{equation}%
Obviously, $\mathcal{N}_{\varepsilon }$ is the square of an equivalent norm
on $\mathcal{H}_{\varepsilon }$, i.e., there is a constant $C>0$,
independent of $\varepsilon $, such that
\begin{equation}
C^{-1}\Vert \bar{\varphi}\Vert _{\mathcal{H}_{\varepsilon }}^{2}\leq
\mathcal{N}_{\varepsilon }(\bar{\varphi})\leq C\Vert \bar{\varphi}\Vert _{%
\mathcal{H}_{\varepsilon }}^{2}.  \label{eq:diff-decomp-2}
\end{equation}%
Following (\ref{eq:diff-decomp-1}) and (\ref{eq:diff-decomp-2}), we have
that, for all $t\geq 0$,
\begin{equation}  \label{eq:to-C2-L}
\Vert \bar{\varphi}(t)\Vert _{\mathcal{H}_{\varepsilon }}^{2}\leq C\Vert
\bar{\varphi}_{0}\Vert _{\mathcal{H}_{\varepsilon }}^{2}e^{-\omega_7 t}.
\end{equation}%
Set $t^{\ast }:=\max \{t_{1},\frac{1}{\omega_7}\ln \left( 4C\right) \}$.
Then, for all $t\geq t^{\ast }$, (\ref{eq:difference-decomposition-L}) holds
with $L_{\varepsilon }=\bar{\varphi}(t^{\ast })$ and
\begin{equation*}
\alpha ^{\ast }=Ce^{-\omega_7 t^{\ast }}<\frac{1}{2}.
\end{equation*}

We now show (\ref{eq:difference-decomposition-K}) holds for some $\Lambda
^{\ast }\geq 0$. First we observe that
\begin{equation}
\begin{aligned}\label{eq:v-bar-identity-2} 2\langle f(v)-f(u),\bar v_{tt}
\rangle_{L^2(\Gamma)} & = \frac{\diff}{\diff t}2\langle f(v)-f(u),\bar v_t
\rangle_{L^2(\Gamma)} - \\ & - 2\langle (f'(v)-f'(u))v_t,\bar v_t
\rangle_{L^2(\Gamma)} + 2\langle f'(u)z_t,\bar v_t \rangle_{L^2(\Gamma)}.
\end{aligned}
\end{equation}%
Next we differentiate the second equation of (\ref{eq:diff-decomp-v}) with
respect to $t,$ multiply the first equation of (\ref{eq:diff-decomp-v}) by $%
2(-\Delta )\bar{v}_{t}$ in $L^{2}(\Omega )$ and insert (\ref%
{eq:v-bar-identity-2}) into the result to produce the differential identity,
which holds for almost all $t\geq 0$,
\begin{equation}
\begin{aligned}\label{eq:v-bar-identity-1} \frac{\diff}{\diff t} & \left\{
\varepsilon\|\bar v_t\|^2_1 + \|\bar v_t\|^2_{L^2(\Gamma)} + \|\Delta\bar
v\|^2 + 2\langle f(u)-f(v),\bar v_t \rangle_{L^2(\Gamma)} \right\} + \\ & +
2\varepsilon\|\bar v_{tt}\|^2_{L^2(\Gamma)} + 2\|\bar v_t\|^2_1 \\ & =
2\langle (f'(v)-f'(u))\nabla v,\nabla\bar v_t \rangle - 2\langle f'(u)\nabla
z,\nabla\bar v_t \rangle + \\ & + 2\langle f(v)-f(u), \bar v_{t}
\rangle_{L^2(\Gamma)} - 2\langle (f'(v)-f'(u))v_t,\bar v_t
\rangle_{L^2(\Gamma)} + 2\langle f'(u)z_t,\bar v_t \rangle_{L^2(\Gamma)}.
\end{aligned}
\end{equation}%
Recall that $z:=u-v$ denotes the difference of any two weak solutions of (%
\ref{eq:pde-h-1})-(\ref{eq:pde-h-3}) and is estimated in (\ref%
{eq:continuous-dependence}). Arguing, for instance, as in \cite[(6.11)-(6.13)%
]{Frigeri10}, we estimate the products on the right hand side of (\ref%
{eq:v-bar-identity-1}), for all $t\in (0,t^{\ast })$, using (\ref%
{eq:assumption-f-1}), Lemma \ref{t:to-H1}, and the embedding $H^{2}(\Omega
)\hookrightarrow C^{0}(\overline{\Omega })$, as follows:
\begin{equation}
\begin{aligned} 2|\langle (f'(u)-f'(v))\nabla v,\nabla\bar v_t \rangle| &
\leq C \left( 1+\|u\|_1 + \|v\|_1 \right)\|z\|_1 \|v\|_2 \|\nabla \bar v_t\|
\\ & \leq C_\varepsilon(t^*) \|\bar{\zeta}_0\|^2_{\mathcal{H}_\varepsilon} +
\frac{1}{4}\|\nabla \bar v_t\|^2, \label{eq:v-bar-estimate-1} \end{aligned}
\end{equation}%
\begin{equation}
\begin{aligned} 2|\langle f'(u)\nabla z,\nabla\bar v_t \rangle| & \leq
C\left( 1+\|u\|^2_2 \right) \|\nabla z\| \|\nabla \bar v_t\| \\ & \leq
C_\varepsilon(t^*) \|\bar{\zeta}_0\|^2_{\mathcal{H}_\varepsilon} +
\frac{1}{4}\|\nabla \bar v_t\|^2, \label{eq:v-bar-estimate-2}\end{aligned}
\end{equation}%
\begin{equation}
\begin{aligned} 2 |\langle f(u)-f(v),\bar v_t \rangle_{L^2(\Gamma)}| & \leq
C \|z\|_1 \|\bar v_t\|_{L^2(\Gamma)} \\ & \leq
C(t^*)\|\bar{\zeta}_0\|^2_{\mathcal{H}_\varepsilon} + \frac{1}{4}\|\bar
v_t\|^2_{L^2(\Gamma)}, \label{eq:v-bar-estimate-3} \end{aligned}
\end{equation}%
\begin{equation}
\begin{aligned} 2|\langle (f'(u)-f'(v))v_t,\bar{v}_t \rangle_{L^2(\Gamma)}|
& \leq C\left( 1+ \|u\|_{C^0(\Gamma)} + \|v\|_{C^0(\Gamma)}
\right)\|z\|_{C^0(\Gamma)}\|v_t\| _{L^2(\Gamma)}\|\bar v_t\|_{L^2(\Gamma)}
\\ & \leq C(t^*)\|\bar{\zeta}_0\|^2_{\mathcal{H}_\varepsilon} +
\frac{1}{4}\|\bar v_t\|^2_{L^2(\Gamma)} \label{eq:v-bar-estimate-4}
\end{aligned}
\end{equation}%
and
\begin{equation}
\begin{aligned} 2|\langle f'(u)z_t,\bar v_t \rangle_{L^2(\Gamma)}| & \leq
C\left( 1+ \|u\|^2_{C^0(\Gamma)} \right)\|z_t\|_{L^2(\Gamma)}\|\bar
v_t\|_{L^2(\Gamma)} \\ & \leq
C(t^*)\|\bar{\zeta}_0\|^2_{\mathcal{H}_\varepsilon} + \frac{1}{4}\|\bar
v_t\|^2_{L^2(\Gamma)}. \label{eq:v-bar-estimate-5} \end{aligned}
\end{equation}%
We emphasize again that by Lemma \ref{t:to-H1} the constants $%
C=C_{\varepsilon }(t^{\ast })$ in estimates (\ref{eq:v-bar-estimate-1}) and (%
\ref{eq:v-bar-estimate-2}) depend on $\varepsilon >0$. After combining (\ref%
{eq:v-bar-estimate-1})-(\ref{eq:v-bar-estimate-5}) with the identity (\ref%
{eq:v-bar-identity-1}), we arrive at the differential inequality,
\begin{equation}
\frac{\mathrm{d}}{\mathrm{d}t}\left\{ \varepsilon \Vert \bar{v}_{t}\Vert
_{1}^{2}+\Vert \bar{v}_{t}\Vert _{L^{2}(\Gamma )}^{2}+\Vert \Delta \bar{v}%
\Vert ^{2}+2\langle f(u)-f(v),\bar{v}_{t}\rangle _{L^{2}(\Gamma )}\right\}
\leq C_{\varepsilon }(t^{\ast })\Vert \bar{\zeta}_{0}\Vert _{\mathcal{H}%
_{\varepsilon }}^{2}  \label{eq:v-bar-estimate-6}
\end{equation}%
(Recall that by the definition of a strong solution in Definition \ref%
{d:regularity-property}, $\bar{v}_{tt}\in L^{2}(0,\infty ;L^{2}(\Gamma ))$).
Now by integrating (\ref{eq:v-bar-estimate-6}) over $(0,t^{\ast })$ and once
again applying the estimate (\ref{eq:v-bar-estimate-3}), we are left with
the bound%
\begin{equation*}
\varepsilon \Vert \bar{v}_{t}(t^{\ast })\Vert _{1}^{2}+\Vert \Delta \bar{v}%
(t^{\ast })\Vert ^{2}\leq C_{\varepsilon }(t^{\ast })\Vert \bar{\zeta}%
_{0}\Vert _{\mathcal{H}_{\varepsilon }}^{2}.
\end{equation*}%
By standard $H^{2}$-elliptic regularity estimates (see (\ref{eq:wt-bounds-1}%
) and (\ref{eq:wt-bounds-3}) above), we obtain
\begin{equation}
\Vert \bar{\theta}(t^{\ast })\Vert _{\mathcal{D}_{\varepsilon }}\leq
C_{\varepsilon }(t^{\ast })\Vert \bar{\zeta}_{0}\Vert _{\mathcal{H}%
_{\varepsilon }}.  \label{eq:to-C2-R}
\end{equation}%
Inequality (\ref{eq:difference-decomposition-K}) now follows with $%
R_{\varepsilon }=\bar{\theta}(t^{\ast })$ and $\Lambda ^{\ast
}=C_{\varepsilon }(t^{\ast })\geq 0$. This finishes the proof.
\end{proof}

\begin{lemma}
\label{t:to-C3} Condition (C3) holds.
\end{lemma}

\begin{proof}
We proceed exactly as in the proof of Lemma \ref{t:to-H1}, differentiating (%
\ref{eq:pde-h-1})-(\ref{eq:pde-h-3}) with respect to $t$ and letting $%
h=u_{t} $. This time we obtain the bound
\begin{equation*}
\Vert \varphi _{t}(t)\Vert _{\mathcal{H}_{\varepsilon }}\leq Q_{\varepsilon
}(R)
\end{equation*}%
for $\varphi _{t}=(u_{t},u_{tt})$, and some function $Q_{\varepsilon }$,
depending on $\varepsilon >0$, where the size of the initial data now
depends on the norm of $\mathcal{B}_{\varepsilon }^{1}$. Hence, on the
compact interval $[t^{\ast },2t^{\ast }]$, the map $t\mapsto S_{\varepsilon
}(t)\varphi _{0}$ is Lipschitz continuous for each fixed $\varphi _{0}\in
\mathcal{B}_{\varepsilon }^{1}$; i.e., there is a constant $L_{\varepsilon
}=L_{\varepsilon }(t^{\ast })>0$ (which depends on $\varepsilon >0$) such
that%
\begin{equation*}
\Vert S_{\varepsilon }(t_{1})\varphi _{0}-S_{\varepsilon }(t_{2})\varphi
_{0}\Vert _{\mathcal{H}_{\varepsilon }}\leq L_{\varepsilon }(t^{\ast
})|t_{1}-t_{2}|.
\end{equation*}%
Together with the continuous dependence estimate (\ref%
{eq:continuous-dependence}), (C3) follows.
\end{proof}

\begin{remark}
\label{rem_att} According to Proposition \ref{abstract1}, the semiflow $%
S_{\varepsilon }:\mathcal{H}_{\varepsilon }\rightarrow \mathcal{H}%
_{\varepsilon }$ possesses an exponential attractor $\mathcal{M}%
_{\varepsilon }\subset \mathcal{B}_{\varepsilon }^{1}$, which attracts
bounded subsets of $\mathcal{B}_{\varepsilon }^{1}$ exponentially fast (in
the topology of $\mathcal{H}_{\varepsilon }$). In order to show that the
attraction property in Theorem \ref{t:exponential-attractors-h}, (iii) also
holds, we can appeal once more to the transitivity of the exponential
attraction \cite[Theorem 5.1]{FGMZ04} and the result of Theorem \ref%
{t:exponential-attraction-h} (also see Remark \ref{expo_att_set}).
\end{remark}

In contrast to the standard case of Dirichlet boundary conditions, where we
have a complete treatment, due to References \cite{FGMZ04, MPZ07}, the
situation with boundary condition (\ref{eq:pde-h-2}) remains essentially
less clear. The following important questions remain open:

\begin{itemize}
\item Higher-order dissipative estimates which are uniform with respect to $%
\varepsilon >0.$

\item Finite-dimensionality of the exponential attractor $\mathcal{M}%
_{\varepsilon }$ (and global attractor $\mathcal{A}_{\varepsilon }$) which
is uniform in $\varepsilon >0.$

\item Existence of a robust (Holder continuous in $\varepsilon \in \left[ 0,1%
\right] $) family of exponential attractors $\left\{ \mathcal{M}%
_{\varepsilon }\right\} .$
\end{itemize}

\section{Appendix}

To make the paper reasonably self-contained, we include the statement of a
frequently used Gr\"{o}nwall-type inequality \cite[Lemma 5]{Pata&Zelik06}.

\begin{proposition}
\label{GL}Let $\Lambda :\mathbb{R}_{+}\rightarrow \mathbb{R}_{+}$ be an
absolutely continuous function satisfying
\begin{equation*}
\frac{\mathrm{d}}{\mathrm{d}t}\Lambda (t)+2\eta \Lambda (t)\leq h(t)\Lambda
(t)+k,
\end{equation*}%
where $\eta >0$, $k\geq 0$ and $\int_{s}^{t}h(\tau )\mathrm{d}\tau \leq \eta
(t-s)+m$, for all $t\geq s\geq 0$ and some $m\geq 0$. Then, for all $t\geq 0$%
,
\begin{equation*}
\Lambda (t)\leq \Lambda (0)e^{m}e^{-\eta t}+\frac{ke^{m}}{\eta }.
\end{equation*}
\end{proposition}

\end{document}